\documentclass[11pt]{article}
\usepackage{graphicx}
\usepackage{amsmath,amssymb}
\newtheorem{theorem}{Theorem}

\def\1{{\bf 1}}

\def\be{\begin{equation}}
\def\ee{\end{equation}}

\begin{document}
\noindent {\large \bf Canonical dual theory applied to a Lennard-Jones potential minimization problem}\\

\bigskip

\noindent {\Large Jiapu Zhang}\\

\noindent Centre for Informatics and Applied Optimization, \&\\ 
Graduate School of Sciences, Information Technology and Engineering,\\  
The University of Ballarat, Mount Helen, VIC 3350, Australia.\\
Emails: j.zhang@ballarat.edu.au, jiapu\_zhang@hotmail.com\\
Telephones: 61-3-5327 9809 (office), 61-4 2348 7360 (mobile)\\

\noindent {\bf Abstract} 
The simplified Lennard-Jones (LJ) potential minimization problem is
\begin{equation*}
\mbox{minimize}~~~f(x)=4\sum_{i=1}^N \sum_{j=1,j<i}^N \left(
\frac{1}{\tau_{ij}^6} -\frac{1}{\tau_{ij}^3}
\right)~~~\mbox{subject to}~~~ x\in \mathbb{R}^n,
\end{equation*}
\noindent where $\tau_{ij}=(x_{3i-2}-x_{3j-2})^2
                +(x_{3i-1}-x_{3j-1})^2
                +(x_{3i}  -x_{3j}  )^2$,
$(x_{3i-2},x_{3i-1},x_{3i})$ is the
coordinates of atom $i$ in $\mathbb{R}^3$, $i,j=1,2,\dots,N(\geq 2 \quad \text{integer})$, $n=3N$ and $N$ is the whole number of atoms. The nonconvexity of the objective function and the huge number of local minima, which is growing exponentially with $N$, interest many mathematical optimization experts. In this paper, the canonical dual theory elegantly tackles this problem illuminated by the amyloid fibril molecular model building.\\


\section{Introduction}
\begin{figure}[h!]
\centerline{
\includegraphics[width=4.2in]{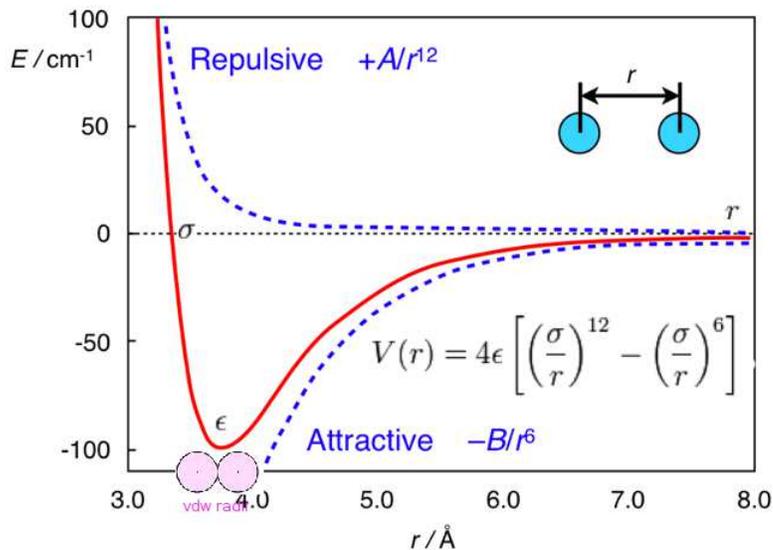}
}
\caption{The Lennard-Jone Potential (formulas (\ref{LJ_r_form}) and (\ref{LJ_AB_form})) (This Figure can be found in website {\tiny http://homepage.mac.com/swain/CMC/DDResources/mol\_interactions/molecular\_interactions.html}}). \label{LJ_potential}
\end{figure}

Neutral atoms are subject to two distinct forces in the limit of large distance and short distance: a dispersion force (i.e. attractive van der Waals (vdw) force) at long ranges, and a repulsion force, the result of overlapping electron orbitals. The Lennard-Jones (L-J) potential represents this behavior ({\small http://en.wikipedia.org/wiki/Lennard-Jones\_potential}, or \cite{locatelli2008} and references therein). The L-J potential is of the form
\begin{equation} \label{LJ_r_form}
V(r)=4\varepsilon \left[ (\frac{\sigma}{r})^{12} - (\frac{\sigma}{r})^6 \right],
\end{equation}
where $r$ is the distance between two atoms, $\varepsilon$ is the depth of the potential well and $\sigma$ is the atom diameter; these parameters can be fitted to reproduce experimental data or deduced from results of accurate quantum chemistry calculations. The $(\frac{\sigma}{r})^{12}$ term describes repulsion and the $(\frac{\sigma}{r})^6$ term describes attraction (Fig. \ref{LJ_potential}). 
In Fig. \ref{LJ_potential} we may see two points: (I) $V(r)=0$ (but the value of $V(r)$ is not the minimal value) when $r=\sigma$ (i.e. the distance between two atoms equals to the sum of {\it atom radii} of the atoms); and (II) when $r=2^{1/6}\sigma$ (i.e. the distance between two atoms equals to the sum of {\it vdw radii} of the atoms), the value of $V(r)$ reaches its minimal value $-\varepsilon$ (i.e. the bottom of the potential well; the force between the atoms is zero at this point). This paper is written based on (II). If we introduce the coordinates of the atoms whose number is denoted by $N$ and let $\varepsilon = \sigma =1$ be the reduced units, the form (\ref{LJ_r_form}) becomes
\begin{equation}\label{LJ_x_form}
f(x)=4\sum_{i=1}^N \sum_{j=1,j<i}^N \left( \frac{1}{\tau_{ij}^6}
-\frac{1}{\tau_{ij}^3} \right),
\end{equation}
\noindent where $\tau_{ij}=(x_{3i-2}-x_{3j-2})^2
                          +(x_{3i-1}-x_{3j-1})^2
                          +(x_{3i}  -x_{3j}  )^2 =||X_i-X_j||^2_2$
and $(x_{3i-2},x_{3i-1},x_{3i})$ is the coordinates of atom $i$, $i,j=1,2,\dots, N (\geq 2)$. The minimization of L-J potential $f(x)$ on $\mathbb{R}^n$ (where $n=3N$) is an optimization problem:
\begin{equation}\label{LJ_f_form}
\min_{s.t. x\in \mathbb{R}^{3N}} f(x).
\end{equation}
This optimization problem interests many optimization experts, for example, Pardalos \cite{pardalossx1994}, Xue \cite{xuemr1992, xue1993, xue1994a, xue1994b, xue2002}, Huang \cite{huang2002b,huang2002a} et al..\\  

For (\ref{LJ_f_form}), when its global optimization solution is reached, the value $r$ in (\ref{LJ_r_form}) should be the sum of two {\it vdw radii} of the two atoms interacted. The three dimensional structure of a molecule with $N$ atoms can be described by specifying the 3-Dimensional coordinate positions $X_1, X_2, \dots, X_N \in \mathbb{R}^3$ of all its atoms. Given bond lengths $r_{ij}$ between a subset $S$ of the atom pairs, the determination of the molecular structure is
\begin{eqnarray}
(\mathcal{P}_0 ) \quad to \quad find \quad &X_1,X_2,\dots ,X_N  \quad s.t. \quad ||X_i-X_j||=r_{ij}, (i,j)\in S,  \label{orginal_problem}
\end{eqnarray}
where $||\cdot ||$ denotes a norm in a real vector space and it is calculated as the Euclidean distance 2-norm in this paper. (\ref{orginal_problem}) can be reformulated as a mathematical global optimization problem (GOP)
\begin{eqnarray}
(\mathcal{P} ) \quad &\min P(X)=\sum_{(i,j)\in S} w_{ij} (||X_i-X_j||^2 -r_{ij}^2 )^2  \label{prime_problem}
\end{eqnarray}
in the terms of finding the global minimum of the function $P(X)$, where $w_{ij}, (i,j)\in S$ are positive weights, $X = (X_1, X_2, \dots, X_N)^T \in \mathbb{R}^{N\times 3}$ \cite{morew1997} and usually $S$ has many fewer than $N^2/2$ elements due to the error in the theoretical or experimental data \cite{grossols2009,zoubs1997}. There may even not exist any solution $X_1, X_2, \dots, X_N$ to satisfy the distance constraints in (\ref{orginal_problem}), for example when data for atoms $i, j, k \in S$ violate the triangle inequality; in this case, we may add a perturbation term $-\epsilon^TX$ to $P(X)$:
\begin{eqnarray}
(\mathcal{P}_{\epsilon} ) \quad &\min P_{\epsilon}(X)=\sum_{(i,j)\in S} w_{ij} (||X_i-X_j||^2 -r_{ij}^2 )^2 -\epsilon^TX, \label{prime_approx_problem}
\end{eqnarray}
where $\epsilon \geq 0$. Thus, the L-J potential optimization problem (\ref{LJ_f_form}) is rewritten into the optimization problem (\ref{prime_approx_problem}).\\

Problem (\ref{prime_approx_problem}) is just the minimization problem of sum of fourth-order polynomials, which can be elegantly solved by the canonical dual theory (CDT) in optimization \cite{gao-book2000,gaorp2010,gaow2012}. We apply the above theory to an amyloid fibril molecular model building problem. The rest of this paper is arranged as follows. In the next section, i.e. Section 2, the CDT will be briefly introduced and its effectiveness will be illuminated by applying the CDT-based optimization approach to the famous double-well problem. In Section 3, the molecular model building works of prion AGAAAAGA amyloid fibrils will be done by the CDT. We find that just very slight refinement is needed by the optimization programs of computational chemistry package Amber 11 to get the optimal molecular models. This implies to us the effectiveness of the CDT to solve our L-J potential minimization problem. Thus, when using the time-consuming and costly X-ray crystallography or NMR spectroscopy we still cannot determine the 3D structure of a protein, we may introduce computational approaches or novel mathematical formulations and physical concepts into molecular biology to study molecular structures. These concluding remarks will be made in the last section, i.e. Section 4.\\

\section{The Canonical dual optimization approach}
We briefly introduce the CDT of \cite{gao-book2000,gaorp2010, gaow2012} specially for solving the following minimization problem of the sum of fourth-order polynomials:
\begin{eqnarray}
&&(\mathcal{P}): \min_x \left\{ P(x) =\sum_{i=1}^m W_i (x)  + \frac{1}{2}x^TQx - x^Tf: x\in \mathbb{R}^n \right\} , \label{primeP}\\
&&where \quad W_i (x)= \frac{1}{2} \alpha_i \left( \frac{1}{2} x^TA_ix +b_i^Tx +c_i \right)^2, A_i=A^T_i, Q=Q^T\in \mathbb{R}^{n\times n}, \nonumber\\ 
&&b_i, f\in \mathbb{R}^n, c_i, \alpha_i \in \mathbb{R}^1, i=1,2,\dots, n, x\in \mathcal{X} \subset \mathbb{R}^n. \nonumber
\end{eqnarray}
The dual problem of $(\mathcal{P})$ is
\begin{eqnarray}
&&(\mathcal{P}^d): \max_{\varsigma} \left\{ P^d(\varsigma ) =\sum_{i=1}^m \left( c_i \varsigma_i -\frac{1}{2} \alpha_i^{-1} \varsigma^2 \right) -\frac{1}{2} F^T(\varsigma ) G^+(\varsigma ) F(\varsigma ): \varsigma \in S_a \right\} ,\label{dualPd}\\
&&where \quad F(\varsigma )=f-\sum_{i=1}^m \varsigma_i b_i, G(\varsigma )=Q+\sum_{i=1}^m \varsigma_i A_i, S_a=\{ \varsigma \in \mathbb{R}^m | F(\varsigma ) \in Col(G(\varsigma )) \}, \nonumber
\end{eqnarray}
$G^+$ denotes the Moore-Penrose generalized inverse of $G$, and $Col(G(\varsigma ))$ is the column space of $G(\varsigma )$. The prime-dual Gao-Strang complementary function of CDT \cite{gao-book2000,gaorp2010, gaow2012} is
\begin{eqnarray}
\Xi (x,\varsigma )=\sum_{i=1}^m \left[ \left( \frac{1}{2} x^TA_ix +b_i^Tx +c_i \right) \varsigma_i -\frac{1}{2} \alpha_i^{-1} \varsigma_i^2 \right] +\frac{1}{2} x^TQx -x^Tf. \label{Gao-StrangComplementary} 
\end{eqnarray}
For $(\mathcal{P})$ and $(\mathcal{P}^d)$ we have the following CDT:
\begin{theorem} \cite{gao-book2000,gaorp2010, gaow2012}
The problem $(\mathcal{P}^d )$ is canonically dual to $(\mathcal{P} )$ in the sense that if $\bar{\varsigma}$ is a critical point of $P^d(\varsigma )$, then $\bar{x} =G^+(\bar{\varsigma} )F(\bar{\varsigma} )$ is a critical point of $P(x)$ on $\mathbb{R}^n$, and $P(\bar{x} )=P^d (\bar{\varsigma} )$. Moreover, if $\bar{\varsigma} \in S_a^+=\{ \varsigma \in S_a | G(\varsigma) \succ 0 \}$, then $\bar{\varsigma}$ is a global maximizer of $P^d(\varsigma )$ over $S_a^+$, $\bar{x}$ is a global minimizer of $P(x)$ on $\mathbb{R}^n$, and
\begin{eqnarray}
P(\bar{x})=\min_{x\in \mathbb{R}^n} P(x)=\Xi (\bar{x}, \bar{\varsigma}) =\max_{\varsigma \in S_a^+} P^d(\varsigma ) =P^d (\bar{\varsigma}).
\end{eqnarray}
\end{theorem}
It is easy to prove that the canonical dual function $P^d(\varsigma )$ is concave on the convex dual feasible space $S_a^+$. Therefore, Theorem 1 shows that the nonconvex primal problem $(\mathcal{P} )$ is equivalent to a concave maximization problem $(\mathcal{P}^d )$ over a convex space $S_a^+$, which can be solved easily by well-developed methods. Over $S_a^-=\{ \varsigma \in S_a | G(\varsigma) \prec 0 \}$ we have the following theorem:
\begin{theorem} \cite{gaow2012}
Suppose that $\bar{\varsigma}$ is a critical point of $(\mathcal{P}^d)$ and the vector $\bar{x}$ is defined by $\bar{x} =G^+(\bar{\varsigma} )F(\bar{\varsigma} )$. If $\bar{\varsigma} \in S_a^-$, then on a neighborhood $\mathcal{X}_o \times \mathcal{S}_o \subset \mathcal{X} \times S_a^-$ of $(\bar{x}, \bar{\varsigma } )$,
we have either
\begin{eqnarray}
P(\bar{x} ) = \min_{x \in \mathcal{X}_o} P(x) =\Xi (\bar{x}, \bar{\varsigma}) =\min_{\varsigma \in \mathcal{S}_o}   P^d(\varsigma ) = P^d(\bar{\varsigma}), \label{equ-dualmin}
\end{eqnarray}
or
\begin{eqnarray}
P(\bar{x} ) = \max_{x \in \mathcal{X}_o} P(x) =\Xi (\bar{x}, \bar{\varsigma}) =\max_{\varsigma \in \mathcal{S}_o}   P^d(\varsigma ) = P^d(\bar{\varsigma}). \label{equ-dualmax}
\end{eqnarray}
\end{theorem}
By the fact that the canonical dual function is a d.c. function (difference of convex functions) on $S_a^-$, the double-min duality (\ref{equ-dualmin}) can be used for finding the biggest local minimizer of $(\mathcal{P} )$ and $(\mathcal{P}^d )$, while the double-max duality (\ref{equ-dualmax}) can be used for finding the biggest local maximizer of $(\mathcal{P} )$ and $(\mathcal{P}^d )$. In physics and material sciences, this pair of biggest local extremal points play important roles in phase transitions.\\

\begin{figure}[h!]
\centerline{
\includegraphics[width=4.2in]{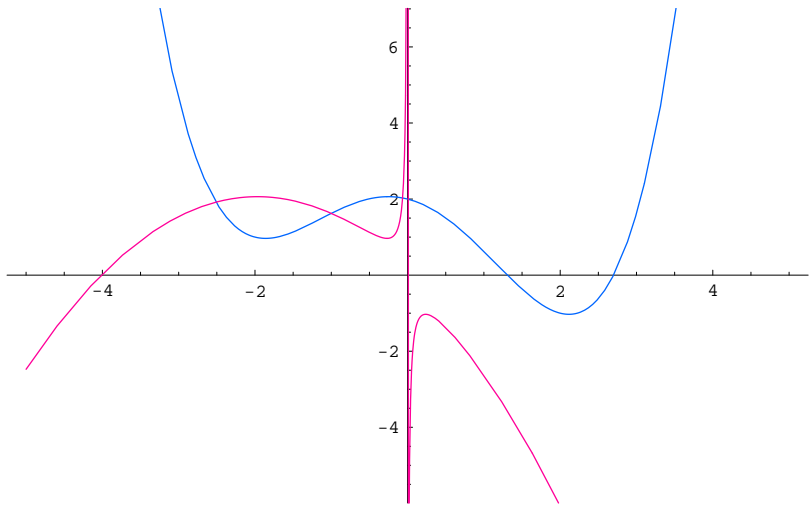}
}
\caption{The prime and dual double-well functions (Prime: blue, Dual: red).}
\label{double_well}
\end{figure}

To illuminate the CDT above-mentioned, we minimize the well-known Double Well potential function \cite{gao-book2000} (blue colored in Fig. \ref{double_well}):
\begin{equation}
P(x)=\frac{1}{2} (\frac{1}{2} x^2- 2)^2 -\frac{1}{2} x. \label{double_well_potential_prime}
\end{equation}
We can easily get $\Xi (x,\varsigma ) =(\frac{1}{2} x^2-2)\varsigma -\frac{1}{2} \varsigma^2 -\frac{1}{2} x$,
\begin{equation}
P^d (\varsigma )=-\frac{1}{8\varsigma }  -\frac{1}{2} \varsigma^2 -2\varsigma \label{double_well_potential_dual}
\end{equation}
(red colored in Fig. \ref{double_well}) and $S_a^+=\{ \varsigma \in \mathbb{R}^1 | \varsigma >0\}$. Let $\Xi (x,\varsigma )'=0$, we get three critical points of $\Xi (x,\varsigma )$: $(\bar{x}^1, \bar{\varsigma}^1)=(2.11491, 0.236417), (\bar{x}^2, \bar{\varsigma}^2)=(-1.86081, -0.268701), (\bar{x}^3, \bar{\varsigma}^3)=(-0.254102,-1.96772)$. By Theorem 1, we know $\bar{x}^1=2.11491$ is the global minimizer of (\ref{double_well_potential_prime}), $\bar{\varsigma}^1=0.236417$ is the global maximizer of (\ref{double_well_potential_dual}) over $S_a^+$, and $P(\bar{\varsigma}^1)=\Xi (\bar{x}^1, \bar{\varsigma}^1)=P^d(\bar{\varsigma}^1)=-1.02951$. By Theorem 2, we know that the local minimizers: $\bar{x}^2=-1.86081, \bar{\varsigma}^2=-0.268701$ (over $S_a^-$), $P(\bar{\varsigma}^2)=\Xi (\bar{x}^2, \bar{\varsigma}^2)=P^d(\bar{\varsigma}^2)=0.9665031$ and the local maximizers: $\bar{x}^3=-0.254102,\bar{\varsigma}^3)=-1.96772$ (over $S_a^-$), $P(\bar{\varsigma}^3)=\Xi (\bar{x}^3, \bar{\varsigma}^3)=P^d(\bar{\varsigma}^3)=2.063$.\\

\section{Applications to a L-J potential optimization problem}
In 2007, Sawaya et al. got a breakthrough finding: the atomic structures of all amyloid fibrils revealed steric zippers, with strong vdw interactions (LJ) between $\beta$-sheets and hydrogen bonds (HBs) to maintain the $\beta$-strands \cite{sawaya2007}. Similarly as (\ref{LJ_r_form}), i.e. the potential energy for the vdw interactions (Fig. \ref{LJ_potential}) between $\beta$-sheets:
\begin{equation} \label{LJ_AB_form}
V_{LJ}(r)=\frac{A}{r^{12}} -\frac{B}{r^6},
\end{equation}
the potential energy for the HBs between the $\beta$-strands has a similar formula 
\begin{equation} \label{HB_r_form}
V_{HB}(r)= \frac{C}{r^{12}} -\frac{D}{r^{10}} ,
\end{equation}
where $A,B,C,D$ are constants given. Thus, the amyloid fibril molecular model building problem is reduced to well solve the optimization problem (\ref{LJ_f_form}) or (\ref{prime_approx_problem}) (in this paper we apply the CDT introduced in Section 2 to solve (\ref{prime_approx_problem})).\\

In this section, we will use suitable templates 3nvf.pdb, 3nvg.pdb and 3nvh.pdb from the Protein Data Bank (http://www.rcsb.org/) to build some amyloid fibril models.\\

\subsection{3NVF}
\begin{figure}[h!]
\centerline{
\includegraphics[scale=0.8]{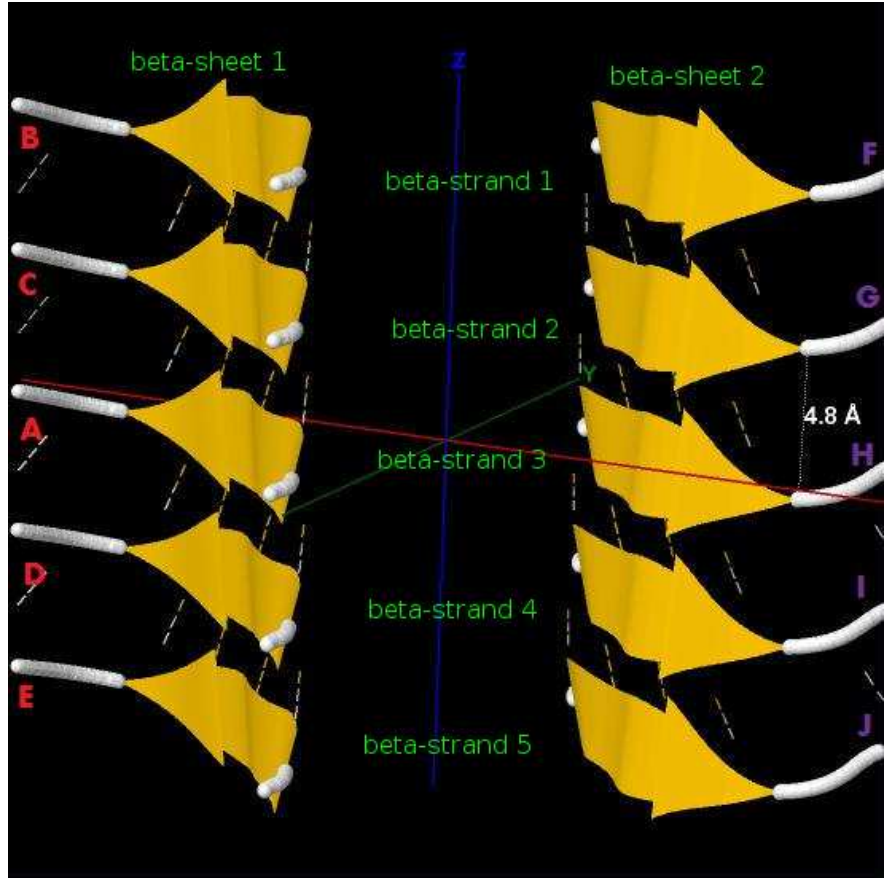}
}
\caption{Protein fibril structure of IIHFGS segment 138--143 from human prion. The purple dashed lines denote the hydrogen bonds. A, B, ..., I, J denote the 10 chains of the fibrils.}
\label{3nvf}
\end{figure}

\begin{figure}[h!]
\centerline{
\includegraphics[scale=0.3]{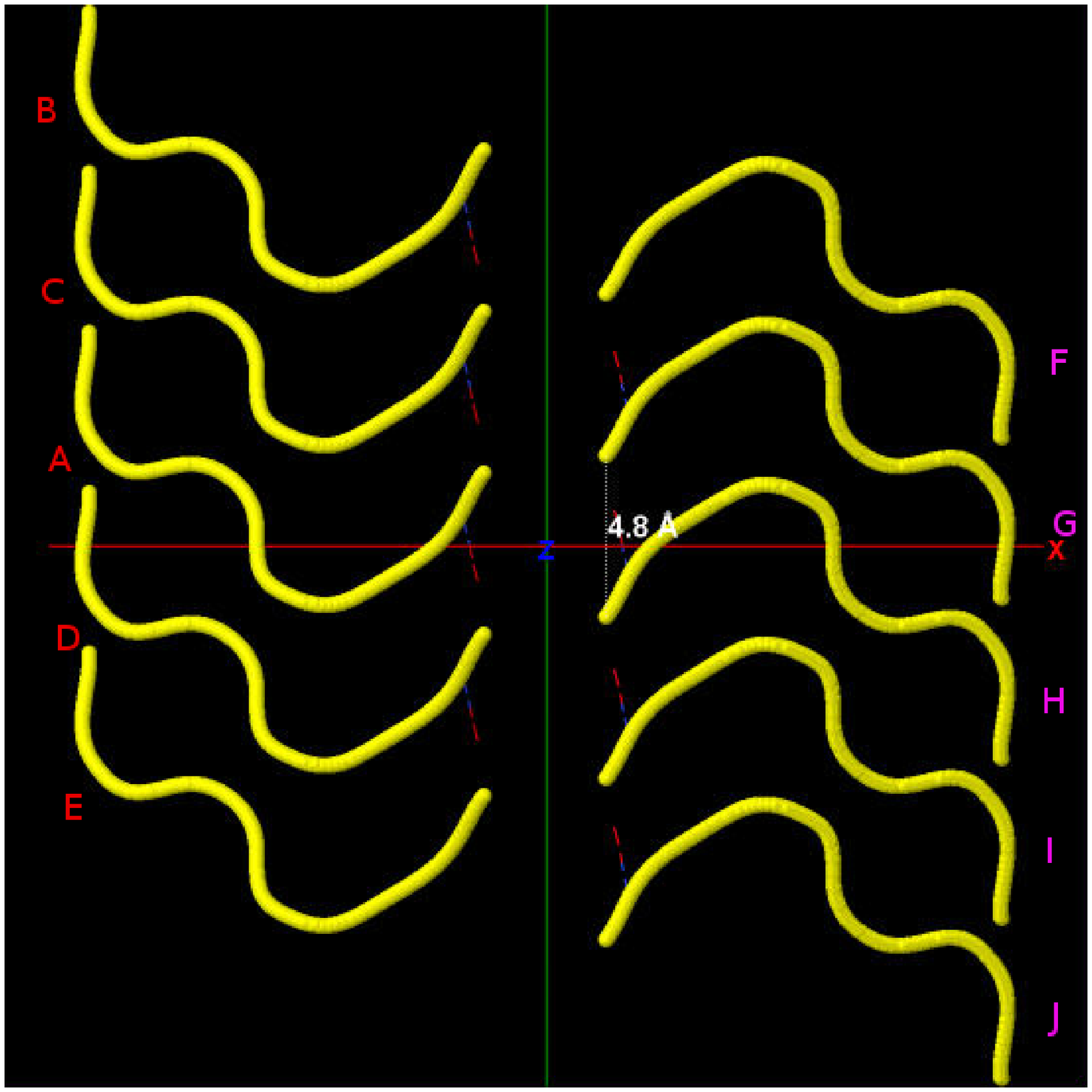}
\includegraphics[scale=0.3]{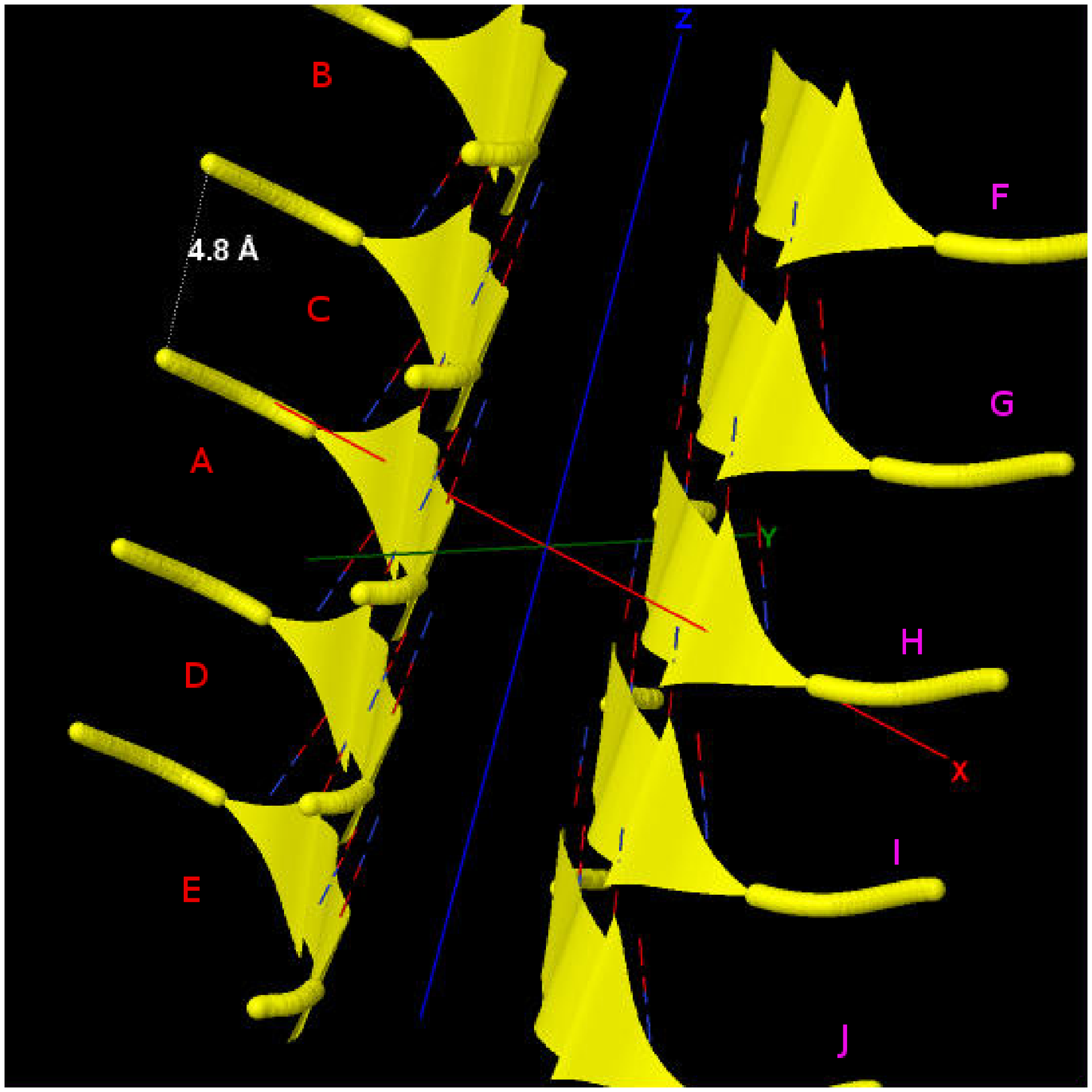}
\includegraphics[scale=0.3]{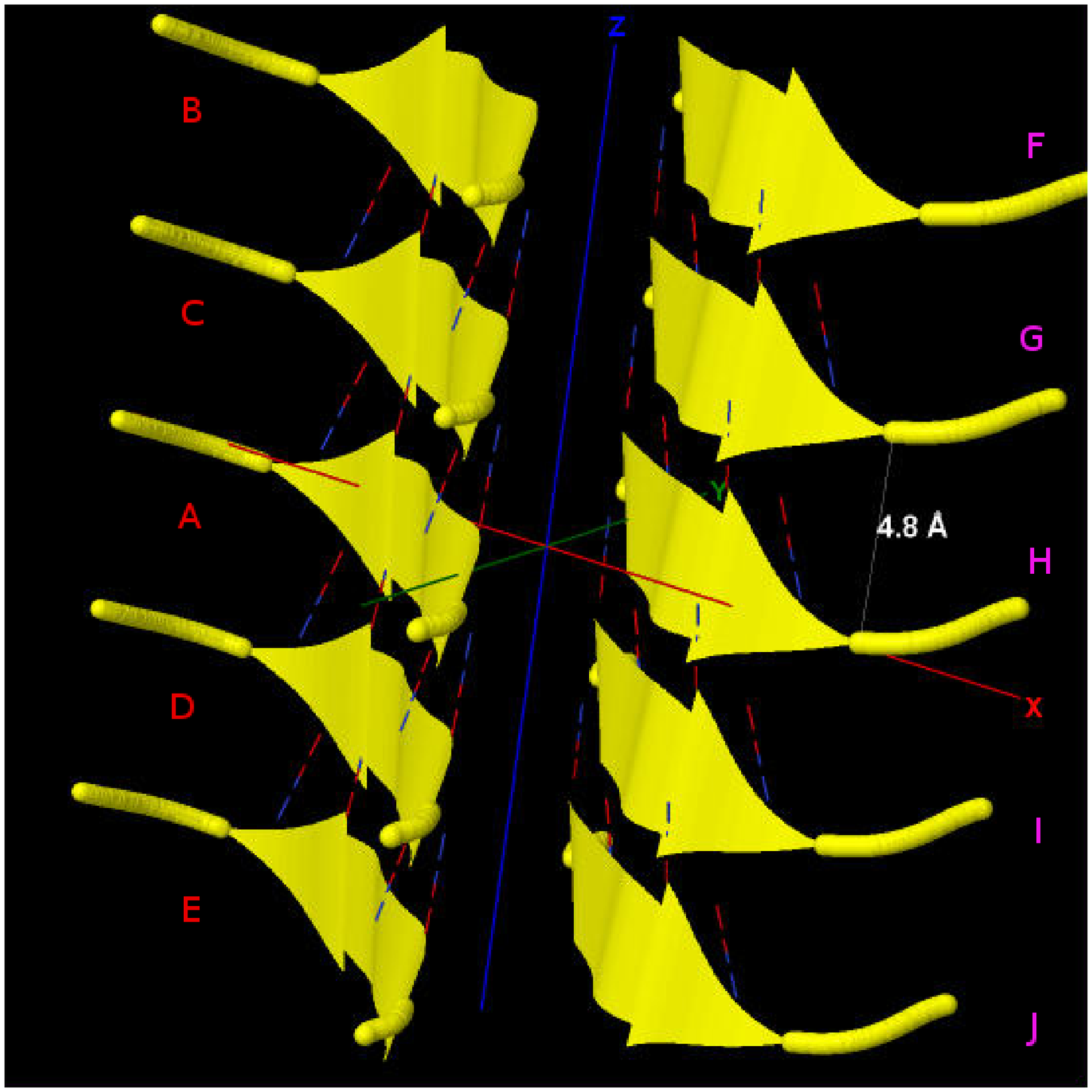}
}
\caption{Protein fibril structure of 3nvf-Models 1-3 (from left to right respectively) for prion AGAAAAGA segment 113-120 . The purple dashed lines denote the hydrogen bonds. A, B, ..., I, J denote the 10 chains of the fibrils.}
\label{3nvf_CDT_models}
\end{figure}

Constructions of the AGAAAAGA amyloid fibril molecular structures of prion 113--120 region are based on the most recently released experimental molecular structures of IIHFGS segment 138--143 from human prion (PDB entry 3NVF released into Protein Data Bank (http://www.rcsb.org) on 2011-03-02) \cite{apostolwsce2011}. The atomic-resolution structure of this peptide is a steric zipper, with strong vdw interactions between $\beta$-sheets and HBs to maintain the $\beta$-strands (Fig. \ref{3nvf}).
In Fig. \ref{3nvf} we see that H chain (i.e. $\beta$-sheet 2) of 3NVF.pdb can be obtained from A chain (i.e. $\beta$-sheet 1) by
\begin{equation}
H = \left( \begin{array}{ccc}
-1   &0   &0\\
 0   &-1  &0\\
 0   &0   &1 
\end{array} \right) A + 
\left( \begin{array}{c}
 27.546\\
 0\\
 0 
\end{array} \right),
\end{equation}
and other chains can be got by
\begin{equation}
C(G)= A (H)+ \left( \begin{array}{c}
0\\
0\\ 
4.8\end{array} \right),
B (F) = A (H) +2\left( \begin{array}{c}
0\\
0\\ 
4.8\end{array} \right), \label{3nvf_cg-bf}
\end{equation}
\begin{equation}
D (I) = A (H) -\left( \begin{array}{c} 
0\\
0\\ 
4.8\end{array} \right),
E (J) = A (H) -2 \left( \begin{array}{c} 
0\\
0\\ 
4.8\end{array} \right). \label{3nvf_di-ej}
\end{equation}
Basing on the template 3NVF.pdb from the Protein Data Bank, three prion AGAAAAGA palindrome amyloid fibril models –- an AGAAAA model (3nvf-Model 1), a GAAAAG model (3nvf-Model 2), and an AAAAGA model (3nvf-Model 3) –- will be successfully constructed in this paper. Chain A of 3nvf-Models 1-3 were got from A Chain of 3NVF.pdb using the mutate module of the free package Swiss-PdbViewer (SPDBV Version 4.01) ({\small http://spdbv.vital-it.ch}). It is pleasant to see that almost all the hydrogen bonds are still kept after the mutations; thus we just need to consider the vdw contacts only. Making mutations for H Chain of 3NVF.pdb, we can get H Chain of 3nvf-Models 1-3. However, we find that the vdw contacts between A Chain and H Chain are too far at this moment. We know that for 3nvf-Model 1 at least the vdw interaction between A.GLY2.CA-H.GLY2.CA, A.ALA4.CB-H.GLY2.CA should be maintained, for 3nvf-Model 2 at least three vdw interactions between A.ALA4.CB-H.ALA2.CB, A.ALA2.CB-H.ALA2.CB, A.ALA2.CB-H.ALA4.CB should be maintained, and for 3nvf-Model 3 at least three vdw interactions between A.ALA2.CB-H.ALA2.CB, A.ALA2.CB-H.ALA4.CB, A.ALA4.CB-H.ALA2.CB should be maintained. Fixing the coordinates of A.GLY2.CA and A.ALA4.CB (two anchors) ((-10.919,-3.862,-1.487), (6.357,1.461,-1.905)) for 3nvf-Model 1, fixing  the coordinates of A.ALA2.CB and A.ALA4.CB (two anchors) ((11.959,-2.844,-1.977), (6.357,1.461,-1.905)) for 3nvf-Models 2-3, letting $d$ equal to the twice of the vdw radius of Carbon atom (i.e. $d =3.4$ angstroms), and letting the coordinates of H.GLY2.CA of 3nvf-Model 1 (two sensors) and the coordinates of H.ALA2.CB and H.ALA4.CB of 3nvf-Models 2-3 (two sensors) be variables, we may get a simple MDGP with 3/6 variables and its dual with 2/3 variables for 3nvf-Model 1:
\begin{eqnarray*}
P_{\epsilon}(x_1)=&&\frac{1}{2} \left\{ (x_{11} +10.919)^2+(x_{12}+3.862)^2 +(x_{13}+1.487)^2 -3.4^2 \right\}^2+\\ 
                      &&\frac{1}{2} \left\{ (x_{11} - 6.357)^2+(x_{12}-1.461)^2 +(x_{13}+1.905)^2 -3.4^2  \right\}^2-\\
                      &&(0.05x_{11}+0.05x_{12}+0.05x_{13}),
\end{eqnarray*}
\begin{eqnarray*}
P^d_{\epsilon}(\varsigma_1, \varsigma_2) =&&124.7908\varsigma_1 -\frac{1}{2} \varsigma_1^2 
                                            +34.615\varsigma_2 -\frac{1}{2} \varsigma_2^2 -\\
&&\frac{1}{2} \left( \begin{array}{c}
0.05-21.838\varsigma_1 +12.714\varsigma_2\\
0.05-7.724\varsigma_1  +2.922\varsigma_2\\
0.05-2.974\varsigma_1   -3.81\varsigma_2
\end{array} \right)^T
\left( \begin{array}{cccccc}
\frac{1}{2\varsigma_1 +2\varsigma_2}  &0                           &0\\
0                           &\frac{1}{2\varsigma_1 +2\varsigma_2}  &0\\
0                           &0                           &\frac{1}{2\varsigma_1 +2\varsigma_2}
\end{array} \right)\\
&&\left( \begin{array}{c}
0.05-21.838\varsigma_1 +12.714\varsigma_2\\
0.05-7.724\varsigma_1  +2.922\varsigma_2\\
0.05-2.974\varsigma_1   -3.81\varsigma_2
\end{array} \right).
\end{eqnarray*}
We can get a global maximal solution (70.1836,70.1812) for $P^d_{\epsilon}(\varsigma_1, \varsigma_2)$ and its corresponding local maximal solution to $P_{\epsilon}(x_1)$:\\
\centerline{$\bar{x}=(-2.28097, -1.20037, -1.69582).$}
By Theorem 1 we know that $\bar{x}$ is a global minimal solution of $P_{\epsilon}(x_1)$. Thus we get 
\begin{equation}
H = \left( \begin{array}{ccc}
-1   &0   &0\\
 0   &-1  &0\\
 0   &0   &1
\end{array} \right) A + 
\left( \begin{array}{c}
-4.5619\\
-2.4009\\
 0.0004
\end{array} \right)
\end{equation}
for 3nvf-Model 1, whose other chains can be got by (\ref{3nvf_cg-bf})-(\ref{3nvf_di-ej}) (Fig. \ref{3nvf_CDT_models}). For 3nvf-Models 2-3, similarly we may get a simple MDGP with 6 variables and its dual with 3 variables:
\begin{eqnarray*}
P_{\epsilon}(x_1,x_2)=&&\frac{1}{2} \left\{ (x_{11} -11.959)^2+(x_{12}+2.844)^2 +(x_{13}+1.977)^2-3.4^2 \right\}^2 +\\ 
                      &&\frac{1}{2} \left\{ (x_{21} -11.959)^2+(x_{22}+2.844)^2 +(x_{23}+1.977)^2 -3.4^2 \right\}^2 +\\
                      &&\frac{1}{2} \left\{ (x_{11}  -6.357)^2+(x_{12}-1.461)^2 +(x_{13}+1.905)^2 -3.4^2 \right\}^2 -\\
                      &&(0.05x_{11}+0.05x_{12}+0.05x_{13}+0.05x_{21}+0.05x_{22}+0.05x_{23}),
\end{eqnarray*}
\begin{eqnarray*}
&&P^d_{\epsilon}(\varsigma_1, \varsigma_2, \varsigma_3 )=143.4545\varsigma_1 -\frac{1}{2} \varsigma_1^2 
                                                        +143.4545\varsigma_2 -\frac{1}{2} \varsigma_2^2 
                                                        + 34.6150\varsigma_3 -\frac{1}{2} \varsigma_3^2 -\\
&&\frac{1}{2} \left( \begin{array}{c}
0.05+23.9180\varsigma_1 +12.7140\varsigma_3\\
0.05- 5.6880\varsigma_1 + 2.9220\varsigma_3\\
0.05- 3.9540\varsigma_1 - 3.8100\varsigma_3\\
0.05+23.9180\varsigma_2\\
0.05- 5.6880\varsigma_2\\
0.05- 3.9540\varsigma_2
\end{array} \right)^T
\left( \begin{array}{cccccc}
\frac{1}{2\varsigma_1 +2\varsigma_3}  &0                           &0                           &0   &0   &0\\
0                           &\frac{1}{2\varsigma_1 +2\varsigma_3}  &0                           &0   &0   &0\\
0                           &0                           &\frac{1}{2\varsigma_1 +2\varsigma_3}  &0   &0   &0\\
0 &0 &0 &\frac{1}{2\varsigma_2} &0 &0\\
0 &0 &0 &0            &\frac{1}{2\varsigma_2} &0\\
0 &0 &0 &0            &0            &\frac{1}{2\varsigma_2}
\end{array} \right)\\
&&\left( \begin{array}{c}
0.05+23.9180\varsigma_1 +12.7140\varsigma_3\\
0.05- 5.6880\varsigma_1 + 2.9220\varsigma_3\\
0.05- 3.9540\varsigma_1 - 3.8100\varsigma_3\\
0.05+23.9180\varsigma_2\\
0.05- 5.6880\varsigma_2\\
0.05- 3.9540\varsigma_2
\end{array} \right) .
\end{eqnarray*}
We can get a global maximal solution (0.920088,0.0127286,0.921273) for $P^d_{\epsilon}(\varsigma_1, \varsigma_2, \varsigma_3)$ and its corresponding local maximal solution to $P_{\epsilon}(x_1,x_2)$:\\
\centerline{$\bar{x}=(9.16977, -0.676538, -1.9274, 13.9231, -0.879925, -0.0129248).$}
By Theorem 1 we know that $\bar{x}$ is a global minimal solution of $P_{\epsilon}(x_1,x_2)$. Thus we get 
\begin{equation}
H = \left( \begin{array}{ccc}
-1   &0   &0\\
 0   &-1  &0\\
 0   &0   &1
\end{array} \right) A + 
\left( \begin{array}{c}
20.8459\\
-2.1533\\
0.6638
\end{array} \right).
\end{equation}
for 3nvg-Models 2-3, whose other chains can be got by (\ref{3nvf_cg-bf})-(\ref{3nvf_di-ej}) (Fig. \ref{3nvf_CDT_models}).
\\

We find 3nvf-Model 1 has some atoms with bad/close contacts, 3nvf-Model 2 has 6 bad/close contacts, and 3nvf-Model 3 has no bad/close contact. This means it not necessary at all to further refine 3nvf-Model 3. We remove these bad contacts by performing energy minimization using Amber 11 \cite{case2010}. Even if there are no obvious bad contacts, it is still a good idea to run a short energy minimization to relax the structures a bit. We will perform the energy minimization in 2 stages. In the first stage, we'll only minimize the water molecules and hold the protein fixed for 500 steps of steepest descent method and then 500 steps of conjugate gradient method. Our goal is just to remove bad contacts, there is no need to go overboard with minimization. In the second stage we proceed directly to minimizing the entire system as a whole for 1500 steps of steepest descent method and then 1000 steps of conjugate gradient method. RMSD (root mean square deviation) is an indicator for structural changes in a protein. It is used to measure the scalar spatial distance between atoms of the same type (for example the C$_\alpha$ atoms) for two structures in different time. The RMSDs between the last snapshot after the refinement and the snapshot illuminated in Fig. \ref{3nvf_CDT_models} are 2.15796, 1.3089087, 1.045318 angstroms for these three 3nvf-Models respectively. The very small values of RMSD are very good measure of precision of CDT for our model building. This shows us that CDT performs well.\\

\subsection{3NVG}
\begin{figure}[h!]
\centerline{
\includegraphics[scale=0.6]{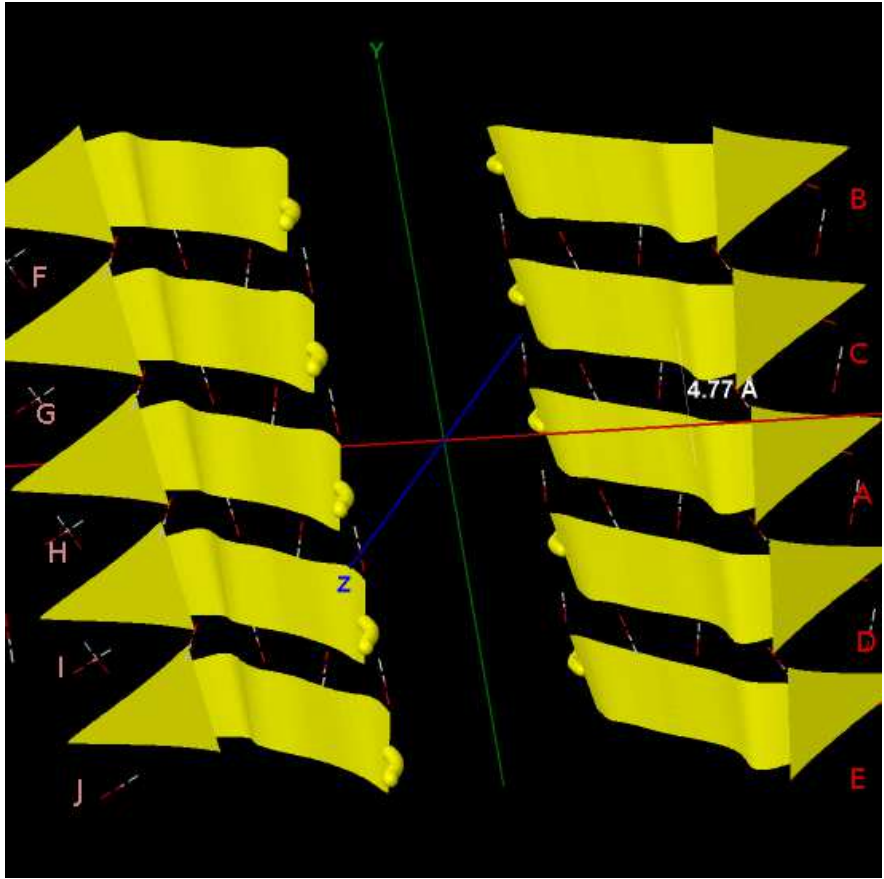}
}
\caption{Protein fibril structure of MIHFGN segment 137--142 from mouse prion. The purple dashed lines denote the hydrogen bonds. A, B, ..., I, J denote the 10 chains of the fibril.}
\label{3nvg}
\end{figure}

\begin{figure}[h!]
\centerline{
\includegraphics[scale=0.3]{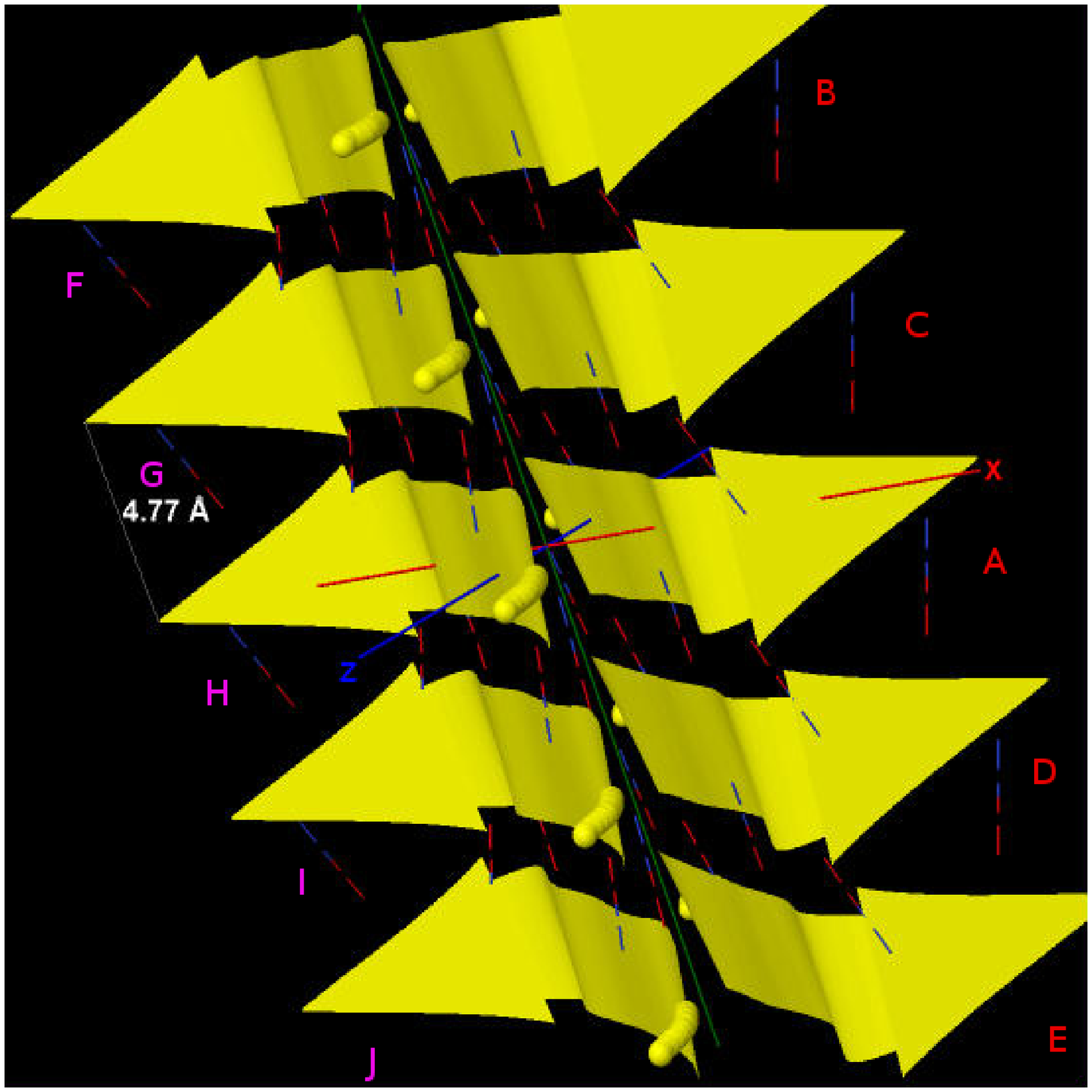}
\includegraphics[scale=0.3]{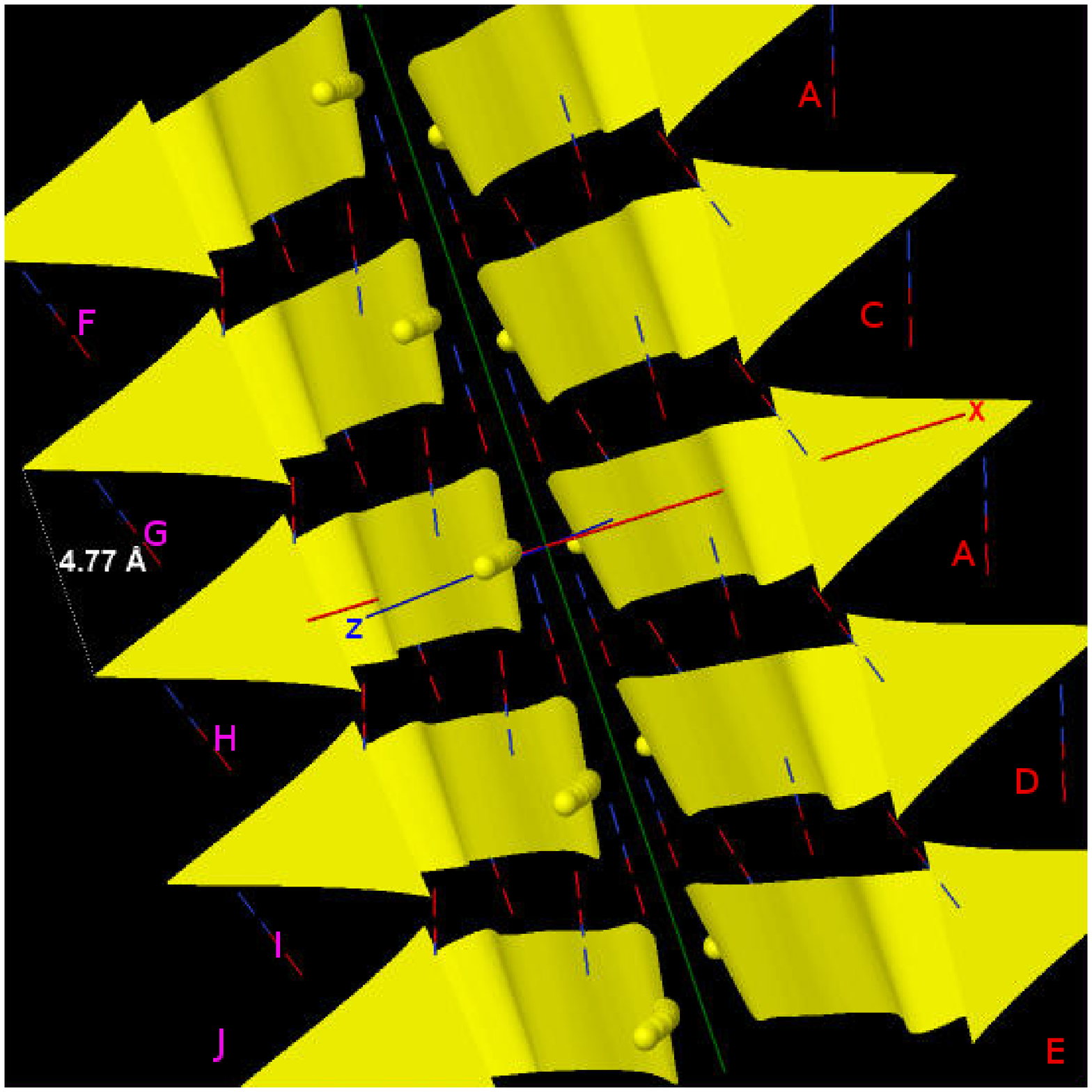}
\includegraphics[scale=0.3]{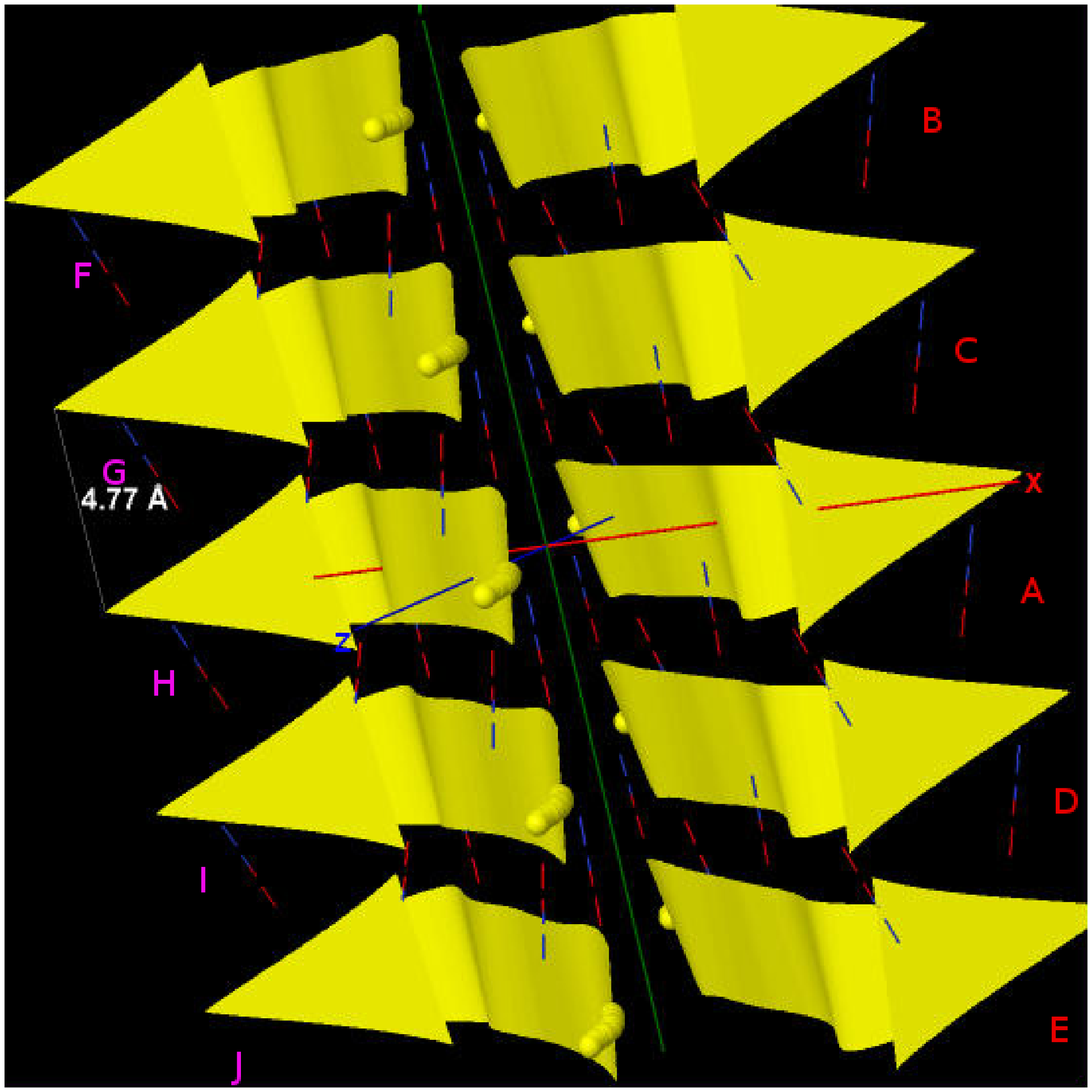}
}
\caption{Protein fibril structure of 3nvg-Models 1-3 (from left to right respectively) for prion AGAAAAGA segment 113-120 . The purple dashed lines denote the hydrogen bonds. A, B, ..., I, J denote the 10 chains of the fibrils.}
\label{3nvg_CDT_models}
\end{figure}

In this subsection, the constructions of the AGAAAAGA amyloid fibril molecular structures of prion 113--120 region are based on the most recently released experimental molecular structures of MIHFGN segment 137--142 from mouse prion (PDB entry 3NVG released into Protein Data Bank (http://www.rcsb.org) on 2011-03-02) \cite{apostolwsce2011}. The atomic-resolution structure of this peptide is a steric zipper, with strong vdw interactions between $\beta$-sheets and HBs to maintain the $\beta$-strands (Fig. \ref{3nvg}).
In Fig. \ref{3nvg} we see that H Chain (i.e. $\beta$-sheet 2) of 3NVG.pdb can be obtained from A Chain (i.e. $\beta$-sheet 1) by
\begin{equation}
H = \left( \begin{array}{ccc}
-1   &0   &0\\
 0   &1   &0\\
 0   &0   &-1 
\end{array} \right) A + 
\left( \begin{array}{c}
-27.28\\
 2.385\\
15.738 
\end{array} \right),
\end{equation}
and other chains can be got by
\begin{equation}
C(G)= A (H)+ \left( \begin{array}{c}
0\\ 
4.77\\
0
\end{array} \right),
B (F) = A (H) +2\left( \begin{array}{c}
0\\ 
4.77\\
0
\end{array} \right), \label{3nvg_cg-bf}
\end{equation}
\begin{equation}
D (I) = A (H) -\left( \begin{array}{c} 
0\\
4.77\\
0
\end{array} \right),
E (J) = A (H) -2 \left( \begin{array}{c} 
0\\ 
4.77\\
0
\end{array} \right). \label{3nvg_di-ej}
\end{equation}
Basing on the template 3NVG.pdb from the Protein Data Bank, three prion AGAAAAGA palindrome amyloid fibril models –- an AGAAAA model (3nvg-Model 1), a GAAAAG model (3nvg-Model 2), and an AAAAGA model (3nvg-Model 3) –- will be successfully constructed in this paper. Chain A of 3nvg-Models 1-3 were got from A Chain of 3NVG.pdb using the mutate module of the free package Swiss-PdbViewer (SPDBV Version 4.01) ({\small http://spdbv.vital-it.ch}). It is pleasant to see that almost all the hydrogen bonds are still kept after the mutations; thus we just need to consider the vdw contacts only. Making mutations for H Chain of 3NVG.pdb, we can get the H Chains of 3nvg-Models 1-3. However, the vdw contacts between A Chain and H Chain are too far at this moment. We know that for 3nvg-Model 1 at least the three vdw interaction between A.GLY2.CA-H.GLY2.CA, A.GLY2.CA-H.ALA4.CB, A.ALA4.CB-H.GLY2.CA should be maintained, for 3nvg-Model 2 at least the three vdw interactions between A.ALA2.CB-H.ALA2.CB, A.ALA2.CB-H.ALA4.CB, A.ALA4.CB-H.ALA2.CB should be maintained, and for 3nvg-Model 3 at least the three vdw interactions between A.ALA2.CB-H.ALA2.CB, A.ALA2.CB-H.ALA4.CB, A.ALA4.CB-H.ALA2.CB should be maintained. Fixing the coordinates of A.GLY2.CA and A.ALA4.CB (two anchors) ((-11.159,-2.241,4.126), (-5.865,-2.618,8.696)) for 3nvg-Model 1, fixing  the coordinates of A.ALA2.CB and A.ALA4.CB (two anchors) ((-12.040,-2.675,5.307), (-5.865,-2.618,8.696)) for 3nvg-Models 1-2, letting $d$ equal to the twice of the vdw radius of Carbon atom (i.e. $d =3.4$ angstroms), and letting the coordinates of H.GLY2.CA and H.ALA4.CB of 3nvg-Model 1 (two sensors) and the coordinates of H.ALA2.CB and H.ALA4.CB of 3nvg-Models 2-3 (two sensors) be variables, we may get a simple MDGP with 6 variables and its dual with 3 variables for 3nvg-Model 1:
\begin{eqnarray*}
P_{\epsilon}(x_1,x_2)=&&\frac{1}{2} \left\{ (x_{11} +11.159)^2+(x_{12}+2.241)^2 +(x_{13}-4.126)^2-3.4^2 \right\}^2+\\ 
                      &&\frac{1}{2} \left\{ (x_{21} +11.159)^2+(x_{22}+2.241)^2 +(x_{23}-4.126)^2 -3.4^2 \right\}^2+\\
                      &&\frac{1}{2} \left\{ (x_{11} + 5.865)^2+(x_{12}+2.618)^2 +(x_{13}-8.696)^2-3.4^2  \right\}^2-\\
                      &&(0.05x_{11}+0.05x_{12}+0.05x_{13}+0.05x_{21}+0.05x_{22}+0.05x_{23}),
\end{eqnarray*}
\begin{eqnarray*}
&&P^d_{\epsilon}(\varsigma_1, \varsigma_2, \varsigma_3) =135.009238\varsigma_1 -\frac{1}{2} \varsigma_1^2 
                                 +135.009238\varsigma_2 -\frac{1}{2} \varsigma_2^2 +105.3125\varsigma_3 -\frac{1}{2} \varsigma_3^2 -\\
&&\frac{1}{2} \left( \begin{array}{c}
0.05-22.318\varsigma_1  -11.73\varsigma_3\\
0.05- 4.482\varsigma_1  -5.236\varsigma_3\\
0.05+ 8.252\varsigma_1 +17.392\varsigma_3\\
0.05-22.318\varsigma_2\\
0.05 -4.482\varsigma_2\\
0.05+ 8.252\varsigma_2
\end{array} \right)^T
\left( \begin{array}{cccccc}
\frac{1}{2\varsigma_1 +2\varsigma_3}  &0                           &0                           &0   &0   &0\\
0                           &\frac{1}{2\varsigma_1 +2\varsigma_3}  &0                           &0   &0   &0\\
0                           &0                           &\frac{1}{2\varsigma_1 +2\varsigma_3}  &0   &0   &0\\
0 &0 &0 &\frac{1}{2\varsigma_2} &0 &0\\
0 &0 &0 &0            &\frac{1}{2\varsigma_2} &0\\
0 &0 &0 &0            &0            &\frac{1}{2\varsigma_2}
\end{array} \right)\\
&&\left( \begin{array}{c}
0.05-22.318\varsigma_1  -11.73\varsigma_3\\
0.05- 4.482\varsigma_1  -5.236\varsigma_3\\
0.05+ 8.252\varsigma_1 +17.392\varsigma_3\\
0.05-22.318\varsigma_2\\
0.05 -4.482\varsigma_2\\
0.05+ 8.252\varsigma_2
\end{array} \right).
\end{eqnarray*}
We can get a global maximal solution (0.708403,0.0127287,0.699001) for $P^d_{\epsilon}(\varsigma_1, \varsigma_2, \varsigma_3)$ and its corresponding local maximal solution to $P_{\epsilon}(x_1,x_2)$:\\
\centerline{$\bar{x}=(-8.51192, -2.41048, 6.4135, -9.19493, -0.276929, 6.09007).$}
By Theorem 1 we know that $\bar{x}$ is a global minimal solution of $P_{\epsilon}(x_1,x_2)$. Thus we get 
\begin{equation}
H = \left( \begin{array}{ccc}
-1   &0   &0\\
 0   &1   &0\\
 0   &0   &-1
\end{array} \right) A + 
\left( \begin{array}{c}
-18.133923\\
0.6673703\\
11.955023 
\end{array} \right)
\end{equation}
for 3nvg-Model 1, whose other chains can be got by (\ref{3nvg_cg-bf})-(\ref{3nvg_di-ej}) (Fig. \ref{3nvg_CDT_models}). For 3nvg-Models 2-3, similarly we may get a simple MDGP with 6 variables and its dual with 3 variables:
\begin{eqnarray*}
P_{\epsilon}(x_1,x_2)=&&\frac{1}{2} \left\{ (x_{11} +12.040)^2+(x_{12}+2.675)^2 +(x_{13}-5.307)^2-3.4^2 \right\}^2 +\\ 
                      &&\frac{1}{2} \left\{ (x_{21} +12.040)^2+(x_{22}+2.675)^2 +(x_{23}-5.307)^2 -3.4^2 \right\}^2 +\\
                      &&\frac{1}{2} \left\{ (x_{11} + 5.865)^2+(x_{12}+2.618)^2 +(x_{13}-8.696)^2-3.4^2 \right\}^2 -\\
                      &&(0.05x_{11}+0.05x_{12}+0.05x_{13}+0.05x_{21}+0.05x_{22}+0.05x_{23}),
\end{eqnarray*}
\begin{eqnarray*}                      
&&P^d_{\epsilon}(\varsigma_1, \varsigma_2, \varsigma_3 )= 168.721474\varsigma_1 -\frac{1}{2} \varsigma_1^2 
                                                        +168.721474\varsigma_2 -\frac{1}{2} \varsigma_2^2 
                                                        +105.312565\varsigma_3 -\frac{1}{2} \varsigma_3^2 -\\
&&\frac{1}{2} \left( \begin{array}{c}
0.05-24.080\varsigma_1  -11.73\varsigma_3\\
0.05-  5.35\varsigma_1  -5.236\varsigma_3\\
0.05+10.614\varsigma_1 +17.392\varsigma_3\\
0.05-24.080\varsigma_2\\
0.05  -5.35\varsigma_2\\
0.05+10.614\varsigma_2
\end{array} \right)^T
\left( \begin{array}{cccccc}
\frac{1}{2\varsigma_1 +2\varsigma_3}  &0                           &0                           &0   &0   &0\\
0                           &\frac{1}{2\varsigma_1 +2\varsigma_3}  &0                           &0   &0   &0\\
0                           &0                           &\frac{1}{2\varsigma_1 +2\varsigma_3}  &0   &0   &0\\
0 &0 &0 &\frac{1}{2\varsigma_2} &0 &0\\
0 &0 &0 &0            &\frac{1}{2\varsigma_2} &0\\
0 &0 &0 &0            &0            &\frac{1}{2\varsigma_2}
\end{array} \right)\\
&&\left( \begin{array}{c}
0.05-24.080\varsigma_1  -11.73\varsigma_3\\
0.05-  5.35\varsigma_1  -5.236\varsigma_3\\
0.05+10.614\varsigma_1 +17.392\varsigma_3\\
0.05-24.080\varsigma_2\\
0.05  -5.35\varsigma_2\\
0.05+10.614\varsigma_2
\end{array} \right) .
\end{eqnarray*}
We can get a global maximal solution (0.849735,0.0127287,0.84036) for $P^d_{\epsilon}(\varsigma_1, \varsigma_2, \varsigma_3)$ and its corresponding local maximal solution to $P_{\epsilon}(x_1,x_2)$:\\
\centerline{$\bar{x}=(-8.95484, -2.63187, 7.00689, -10.0759, -0.710929, 7.27107).$}
By Theorem 1 we know that $\bar{x}$ is a global minimal solution of $P_{\epsilon}(x_1,x_2)$. Thus we get 
\begin{equation}
H = \left( \begin{array}{ccc}
-1   &0   &0\\
 0   &1   &0\\
 0   &0   &-1
\end{array} \right) A + 
\left( \begin{array}{c}
-19.3102\\
0.6644\\
13.5316 
\end{array} \right).
\end{equation}
for 3nvg-Models 2-3, whose other chains can be got by (\ref{3nvg_cg-bf})-(\ref{3nvg_di-ej}) (Fig. \ref{3nvg_CDT_models}).
\\

We did same refinements for 3nvg-Models 1-3 as for 3nvf-Models 1-3. The RMSDs between the last snapshot after the refinement and the snapshot illuminated in Fig. \ref{3nvg_CDT_models} are 1.572438, 1.404648, 1.464767 angstroms for these three 3nvg-Models respectively. The very small values of RMSD  show us that CDT performs well and precisely for 3nvg-Model building.\\

\subsection{3NVH}
\begin{figure}[h!]
\centerline{
\includegraphics[scale=0.6]{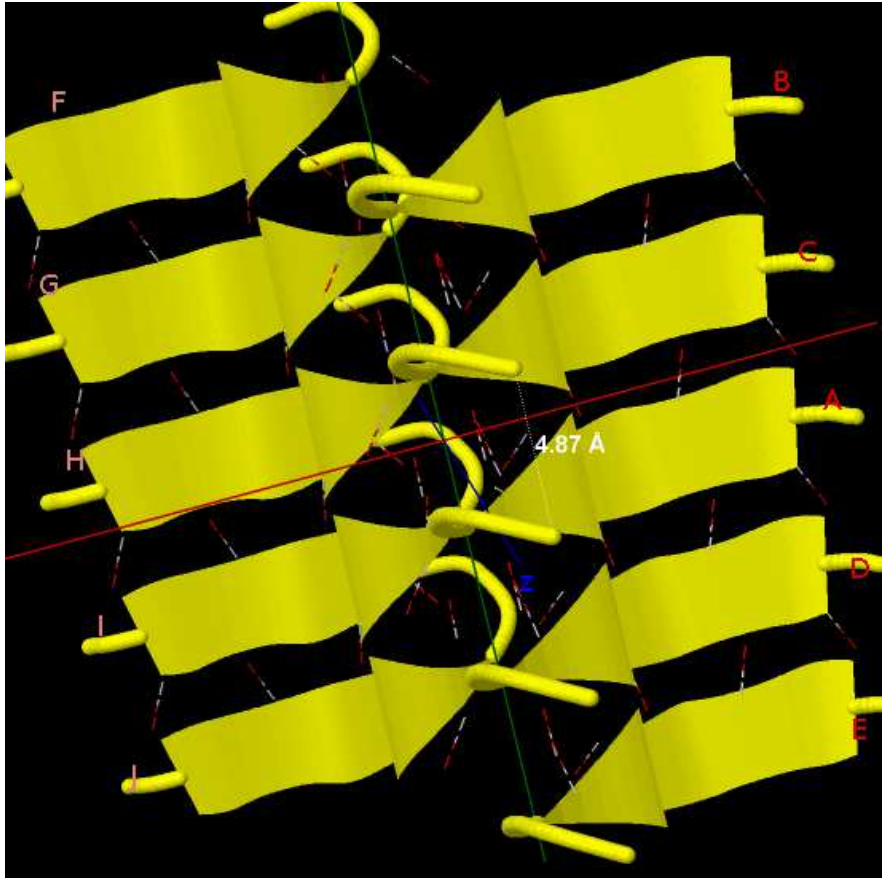}
}
\caption{Protein fibril structure of MIHFGND segment 137--143 from mouse prion. The purple dashed lines denote the hydrogen bonds. A, B, ..., I, J denote the 10 chains of the fibril.}
\label{3nvh}
\end{figure}

\begin{figure}[h!]
\centerline{
\includegraphics[scale=0.3]{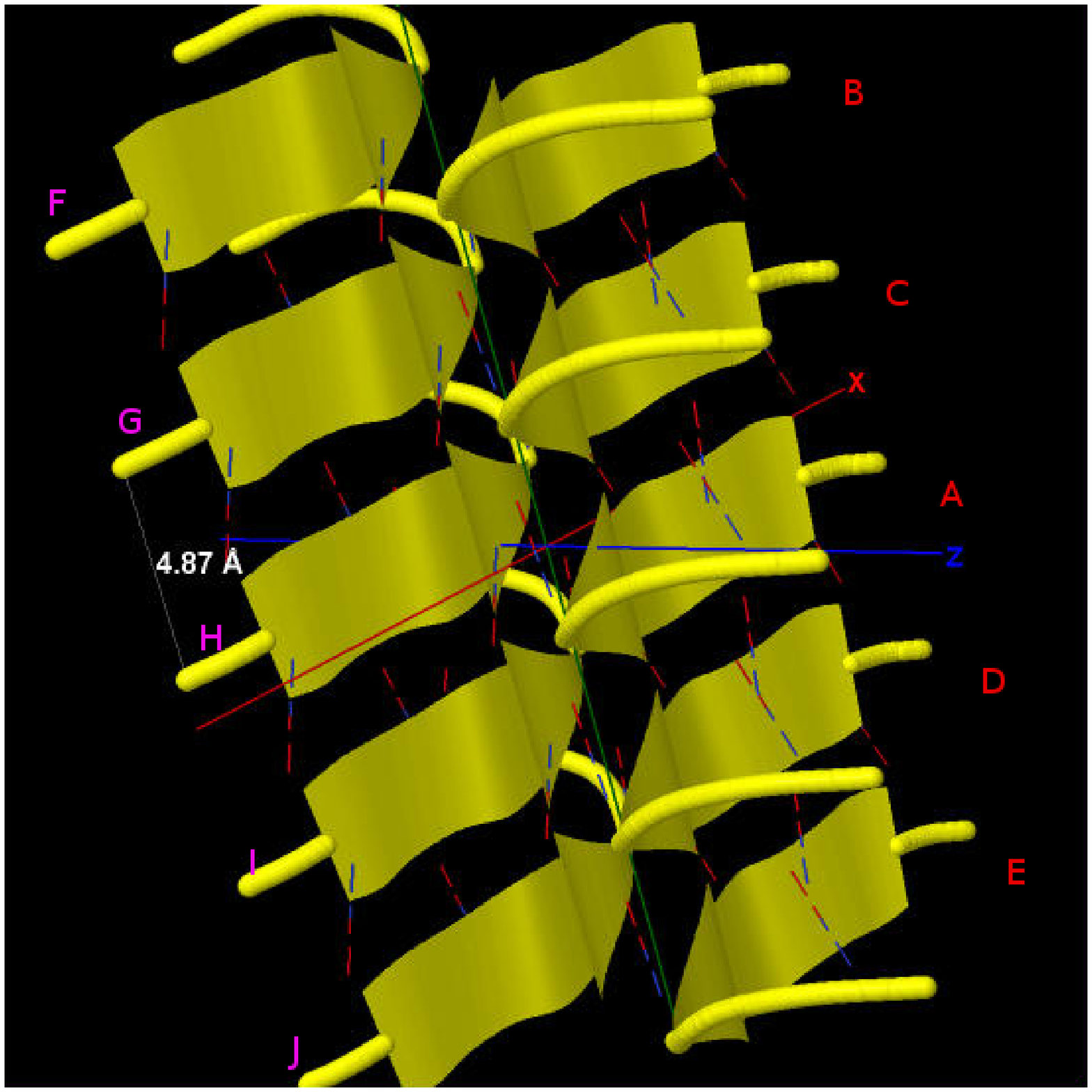}
\includegraphics[scale=0.3]{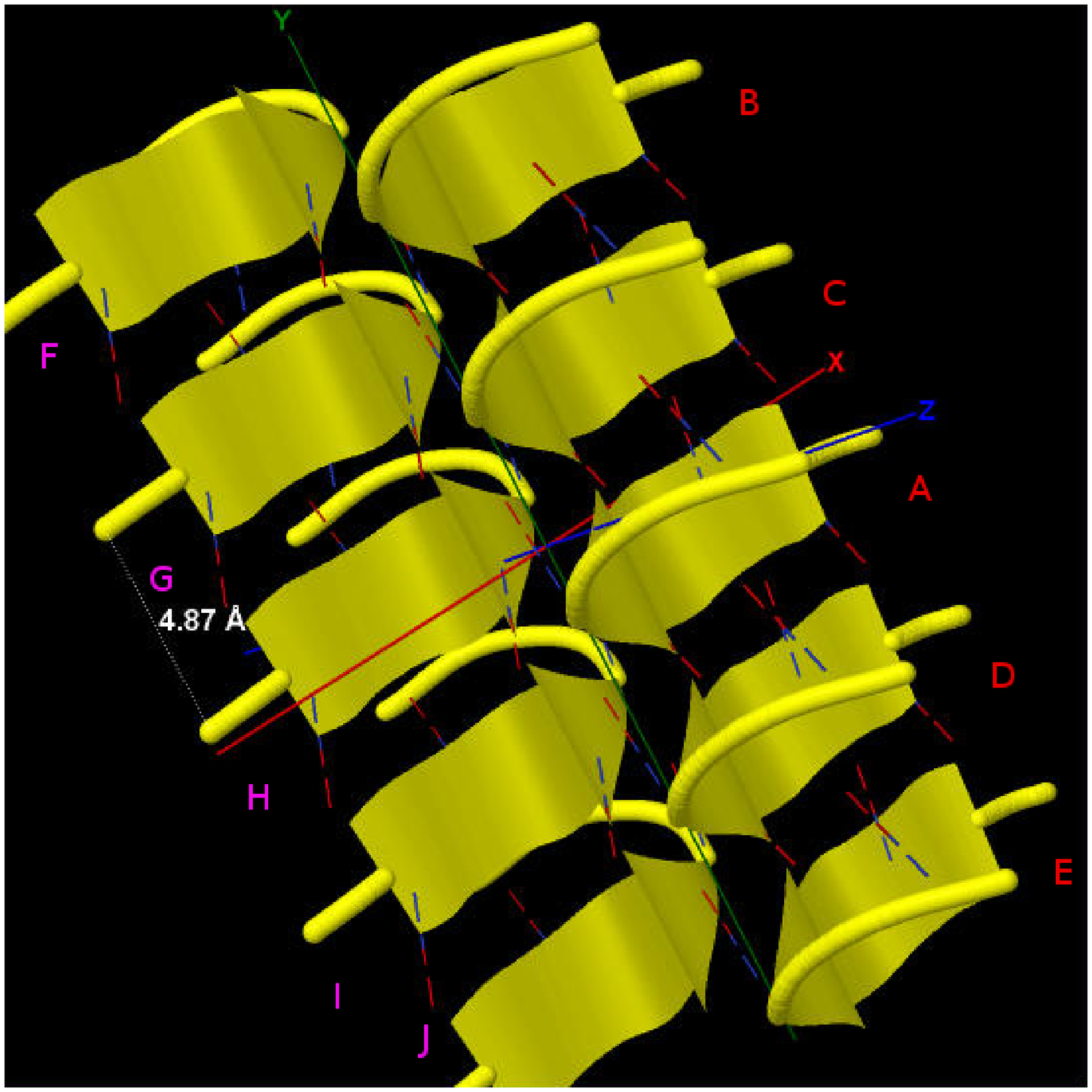}
}
\caption{Protein fibril structure of 3nvh-Models 1-2 (from left to right respectively) for prion AGAAAAGA segment 113-120 . The purple dashed lines denote the hydrogen bonds. A, B, ..., I, J, K, L denote the 12 chains of the fibrils.}
\label{3nvh_CDT_models}
\end{figure}

Similar as the above two subsections, this subsection constructs the AGAAAAGA amyloid fibril molecular structures of prion 113--120 region basing on the most recently released experimental molecular structures of MIHFGND segment 137--143 from mouse prion (PDB entry 3NVH released into Protein Data Bank (http://www.rcsb.org) on 2011-03-02) \cite{apostolwsce2011}. The atomic-resolution structure of this peptide is a steric zipper, with strong vdw interactions between $\beta$-sheets and HBs to maintain the $\beta$-strands (Fig. \ref{3nvh}).
In Fig. \ref{3nvh} we see that H chain (i.e. $\beta$-sheet 2) of 3NVH.pdb can be obtained from A chain (i.e. $\beta$-sheet 1) by
\begin{equation}
H = \left( \begin{array}{ccc}
-1   &0   &0\\
 0   &1   &0\\
 0   &0   &-1
\end{array} \right) A + 
\left( \begin{array}{c}
0\\
2.437\\
-15.553 
\end{array} \right),
\end{equation}
and other chains can be got by
\begin{equation}
C(G)= A (H)+ \left( \begin{array}{c}
0\\ 
4.87\\
0
\end{array} \right), 
B (F) = A (H) +2\left( \begin{array}{c}
0\\ 
4.87\\
0
\end{array} \right), \label{3nvh_cg-bf}
\end{equation}
\begin{equation}
D (I) = A (H) -\left( \begin{array}{c} 
0\\
4.87\\
0
\end{array} \right),
E (J) = A (H) -2 \left( \begin{array}{c} 
0\\ 
4.87\\
0
\end{array} \right). \label{3nvh_di-ej}
\end{equation}
Basing on the template 3NVH.pdb from the Protein Data Bank, three prion AGAAAAGA palindrome amyloid fibril models –- an AGAAAAG model (3nvh-Model 1), a GAAAAGA model (3nvh-Model 2) –- will be successfully constructed in this paper. A chain of 3nvh-Models 1-2 were got from A chain of 3NVH.pdb using the mutate module of the free package Swiss-PdbViewer (SPDBV Version 4.01) ({\small http://spdbv.vital-it.ch}). It is pleasant to see that almost all the hydrogen bonds are still kept after the mutations; thus we just need to consider the vdw contacts only. Making mutations for H chain of 3NVH.pdb, we can get the H chains of 3nvh-Models 1-2. However, the vdw contacts between A chain and H chain are too far at this moment ($\geq$4.25 Angstroms). we may know that for 3nvh-Models 1-2 at least one vdw interaction between A.ALA4.CB-H.ALA4.CB should be maintained. Fixing the coordinates of A.ALA4.CB (the anchor) ((1.731,-1.514,-7.980)), letting $d$ equal to the twice of the vdw radius of Carbon atom (i.e. $d =3.4$ angstroms), and letting the coordinate of H.ALA4.CB (one sensor) be variables, we may get a simple MDGP with 3 variables and its dual with 1 variable:
\begin{eqnarray*}
P_{\epsilon}(x_1)=  &&\frac{1}{2} \left\{ (x_{11} -1.731)^2+(x_{12}+1.514)^2 +(x_{13}+7.980)^2 -3.4^2 \right\}^2\\
                    && -0.05x_{11}-0.05x_{12}-0.05x_{13},\\
P^d_{\epsilon}(\varsigma_1) = &&57.409\varsigma_1 -\frac{1}{2} \varsigma_1^2\\
                    &&\frac{(0.05 + 3.462\varsigma_1)^2+(0.05 - 3.028\varsigma_1)^2+(0.05 - 15.96\varsigma_1)^2}{4\varsigma_1} .
\end{eqnarray*}
We can easily get the global maximal solution $ 0.0127287 \in \{ \varsigma \in \mathbb{R}^1 | \varsigma_i >0, i=1\}$ for $P_{\epsilon}^d(\varsigma_1)$. Then, we get its corresponding solution for $P_{\epsilon}(x_1)$:\\
\centerline{$\bar{x}=(3.69507, 0.450071, -6.01593).$}
By Theorem 1 we know that $\bar{x}$ is a global minimal solution of $P_{\epsilon}(x_1)$, i.e. for H.ALA4.CB. Thus we get 
\begin{equation}
H = \left( \begin{array}{ccc}
-1   &0   &0\\
 0   &1   &0\\
 0   &0   &-1
\end{array} \right) A + 
\left( \begin{array}{c}
5.42607\\
1.964071\\
-13.99593 
\end{array} \right)
\end{equation}
for 3nvh-Models 1-2, whose other chains can be got by (\ref{3nvh_cg-bf})-(\ref{3nvh_di-ej}) (Fig. \ref{3nvh_CDT_models}).
\\

We carried on the same refinements for 3nvh-Models 1-2 as for 3nvf-Models 1-3 and 3nvg-Models 1-3. The RMSDs between the last snapshot after the refinement and the snapshot illuminated in Fig. \ref{3nvh_CDT_models} are 1.534417, 1.572836 angstroms for the two 3nvh-Models respectively. The very small values of RMSD again show to us that CDT performs well and precisely for 3nvg-Model building.\\

\subsection{Refined 3nvf-Models 1-3, 3nvg-Models 1-3, 3nvh-Models 1-2}
\begin{figure}[h!]
\centerline{
\includegraphics[scale=0.2]{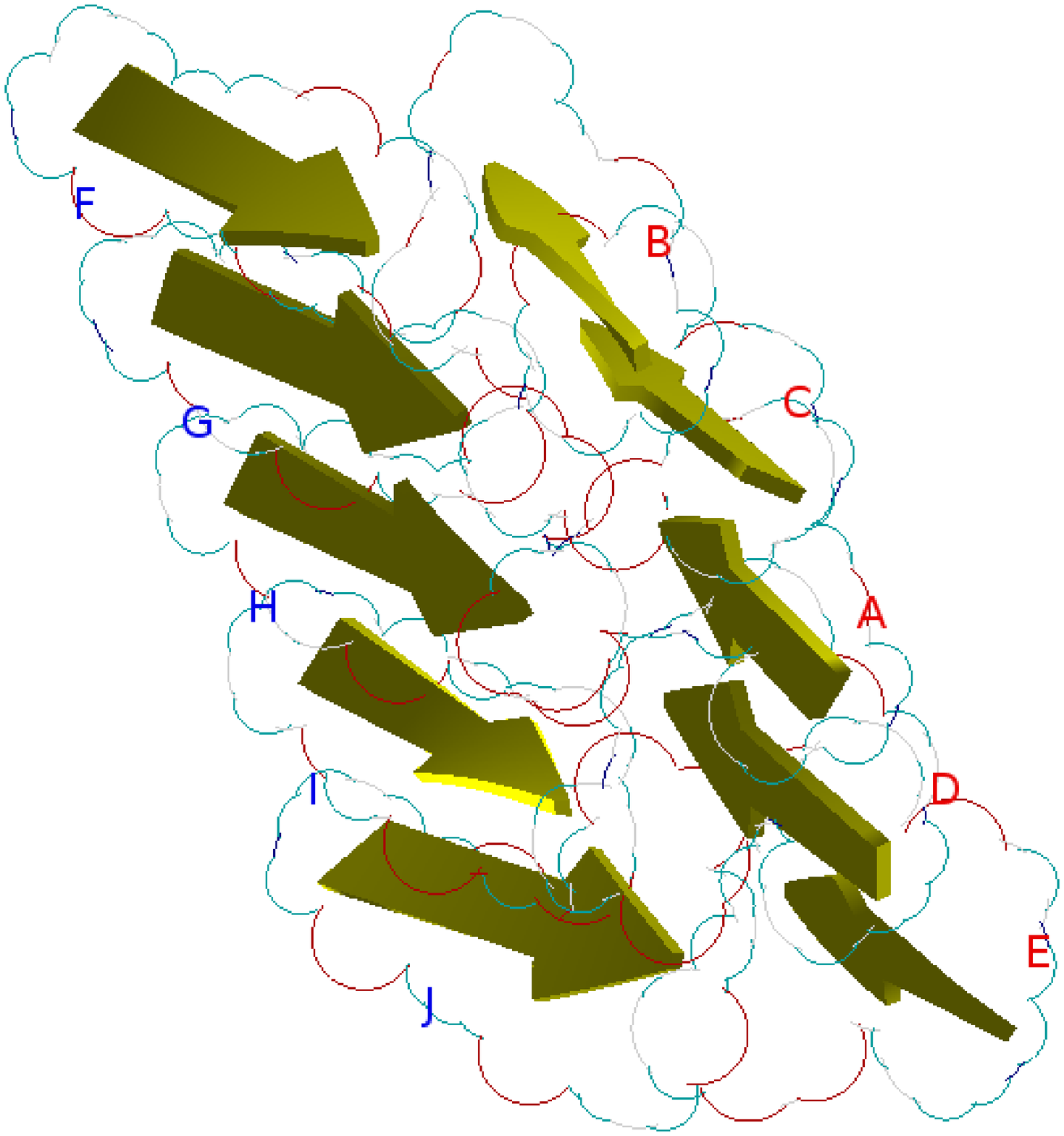}
\includegraphics[scale=0.2]{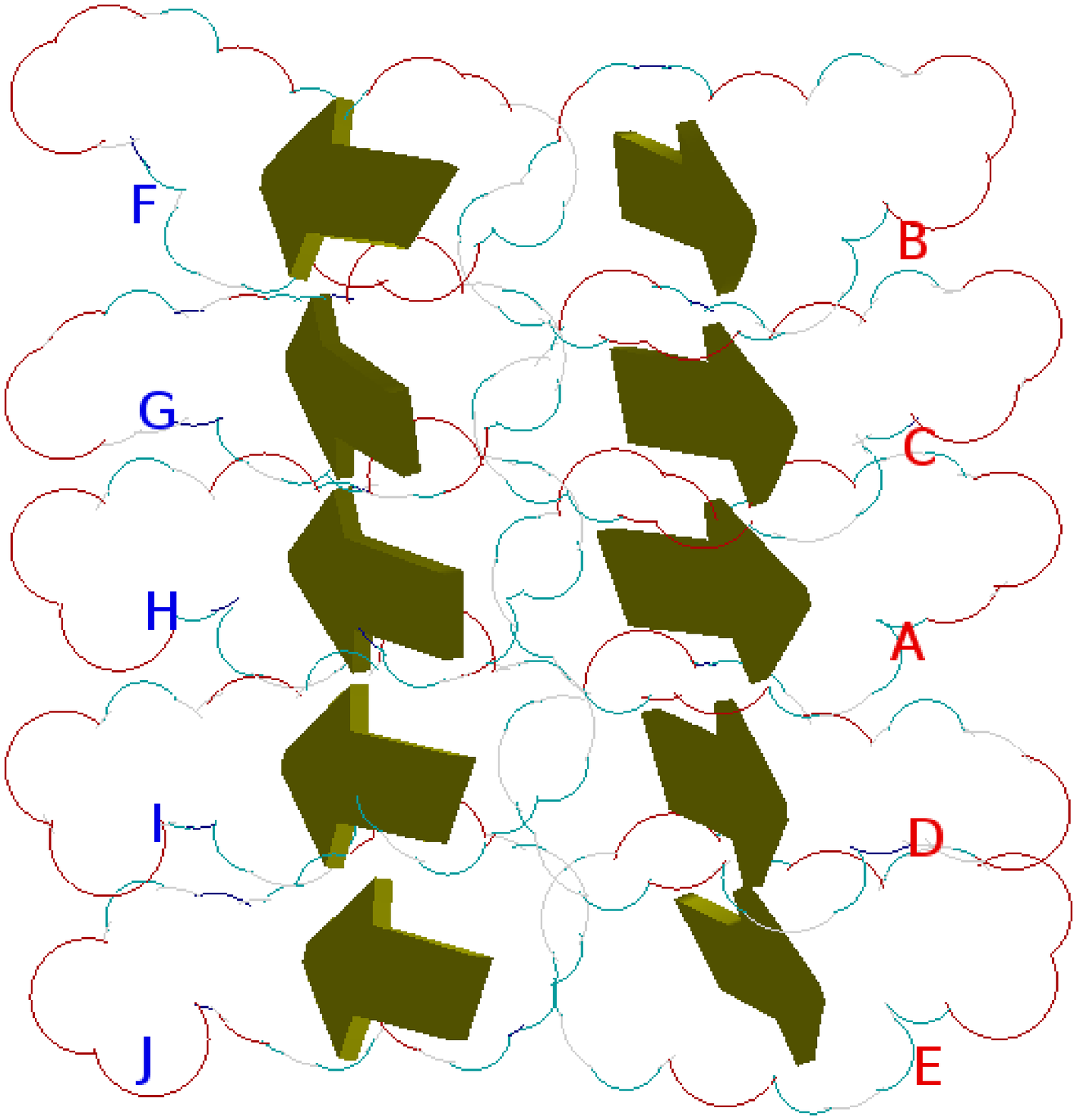}
\includegraphics[scale=0.2]{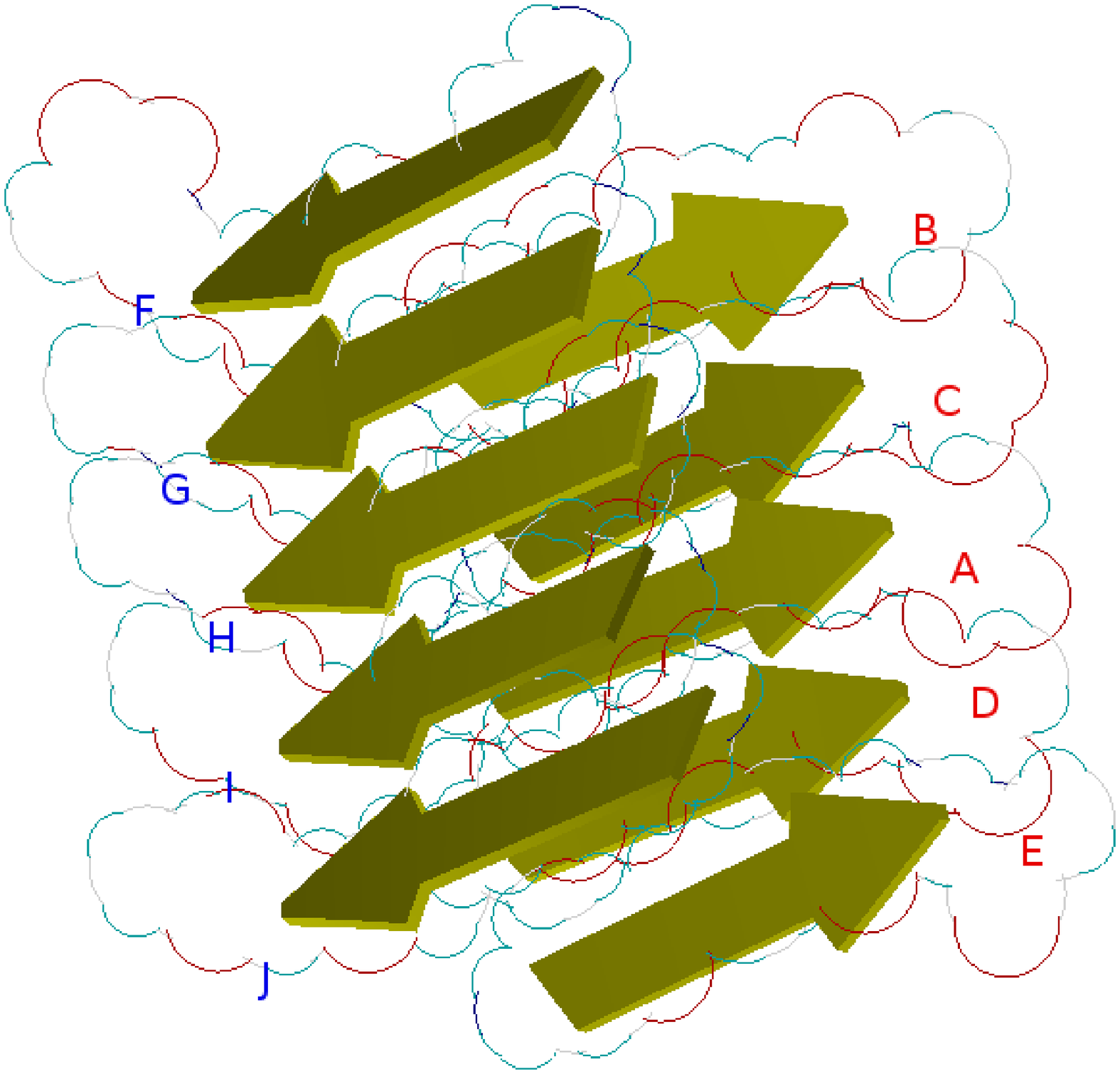}
}
\caption{Perfect 3nvf-Models 1-3 (from left to right respectively) for prion AGAAAAGA segment 113-120 . The purple dashed lines denote the hydrogen bonds. A, B, ..., I, J denote the 10 chains of the fibrils.}
\label{3nvf_CDT_models_min2}
\end{figure}

\begin{figure}[h!]
\centerline{
\includegraphics[scale=0.2]{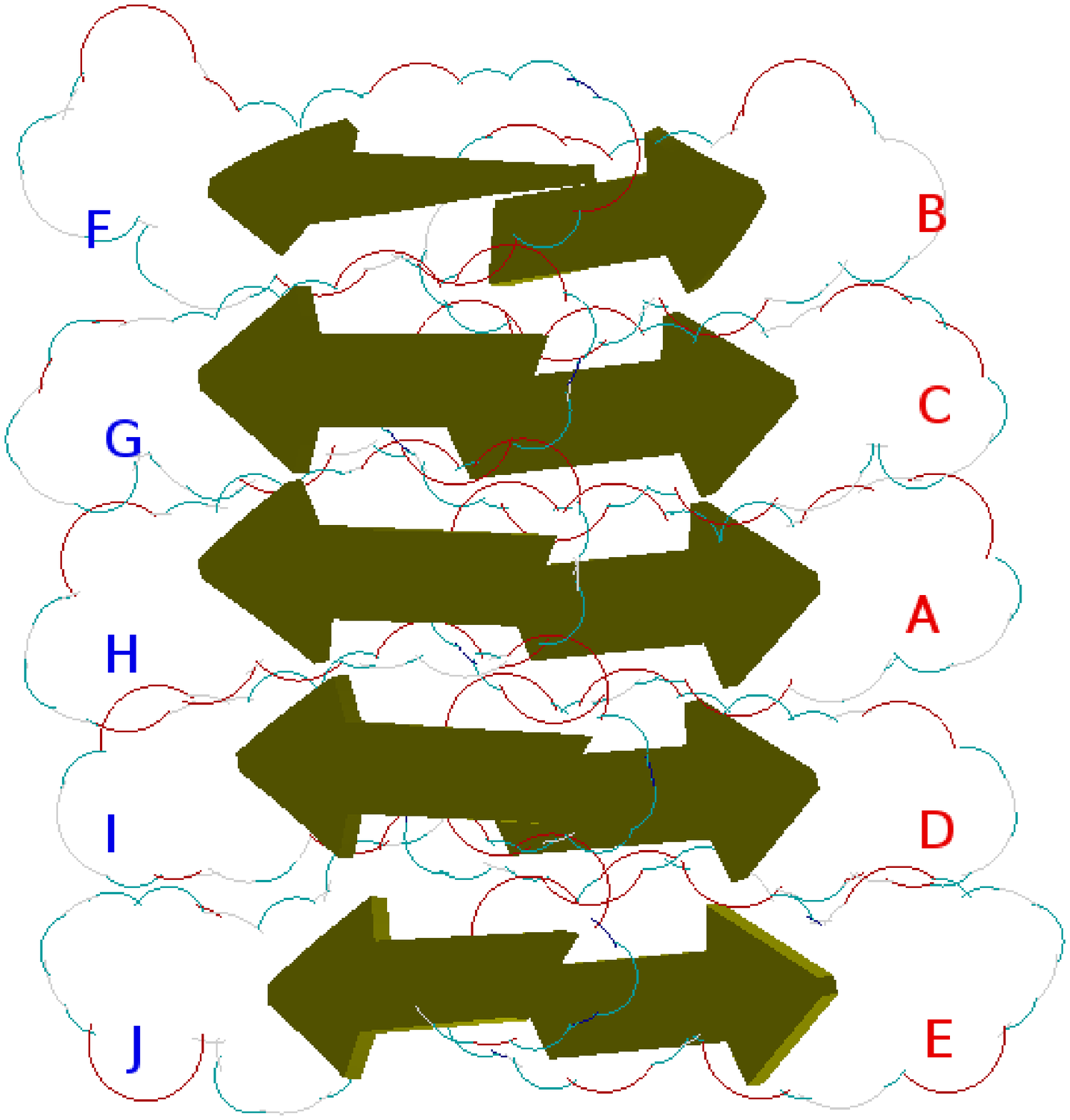}
\includegraphics[scale=0.2]{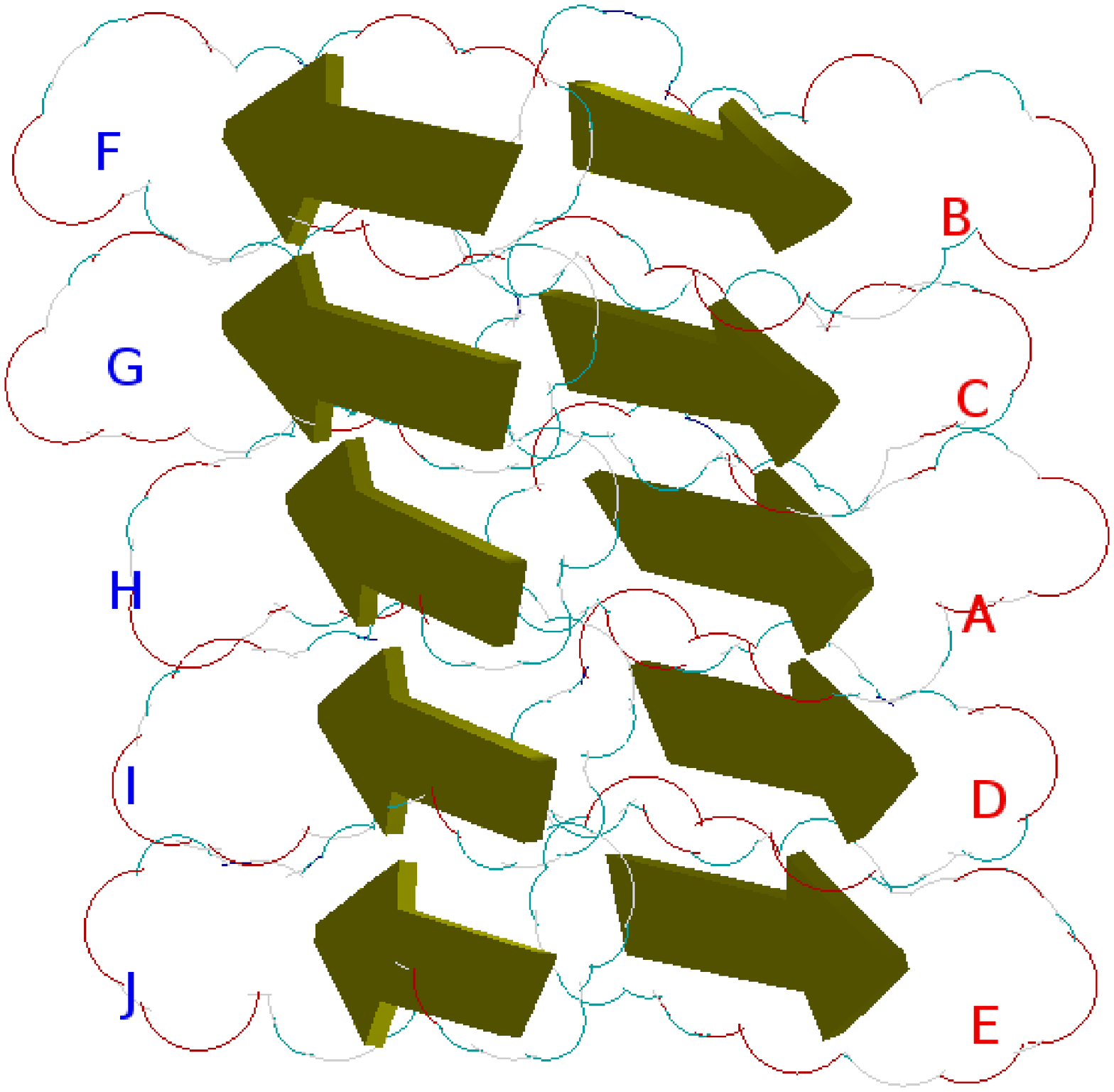}
\includegraphics[scale=0.2]{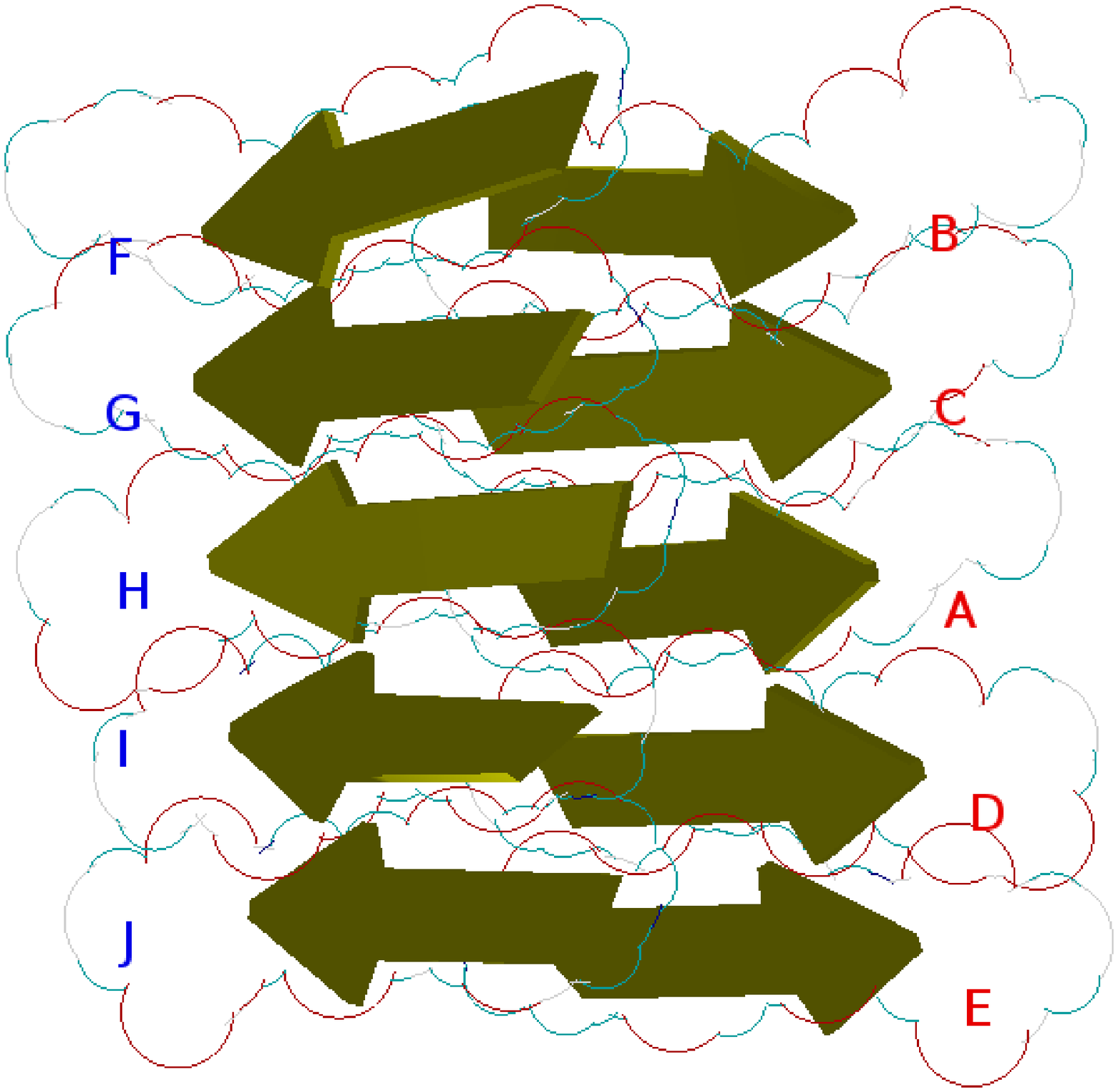}
}
\caption{Perfect 3nvg-Models 1-3 (from left to right respectively) for prion AGAAAAGA segment 113-120 . The purple dashed lines denote the hydrogen bonds. A, B, ..., I, J denote the 10 chains of the fibrils.}
\label{3nvg_CDT_models_min2}
\end{figure}

\begin{figure}[h!]
\centerline{
\includegraphics[scale=0.2]{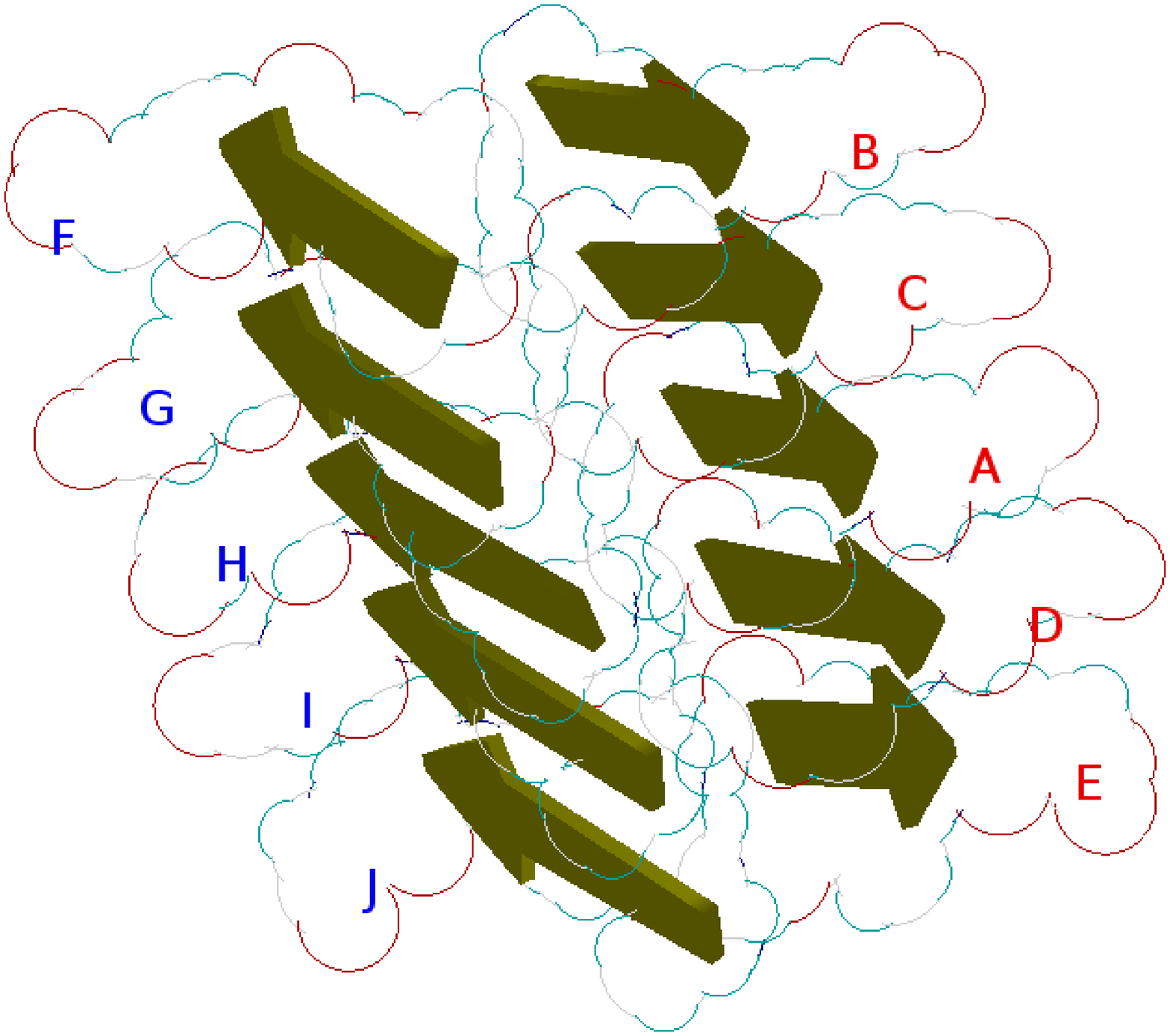}
\includegraphics[scale=0.2]{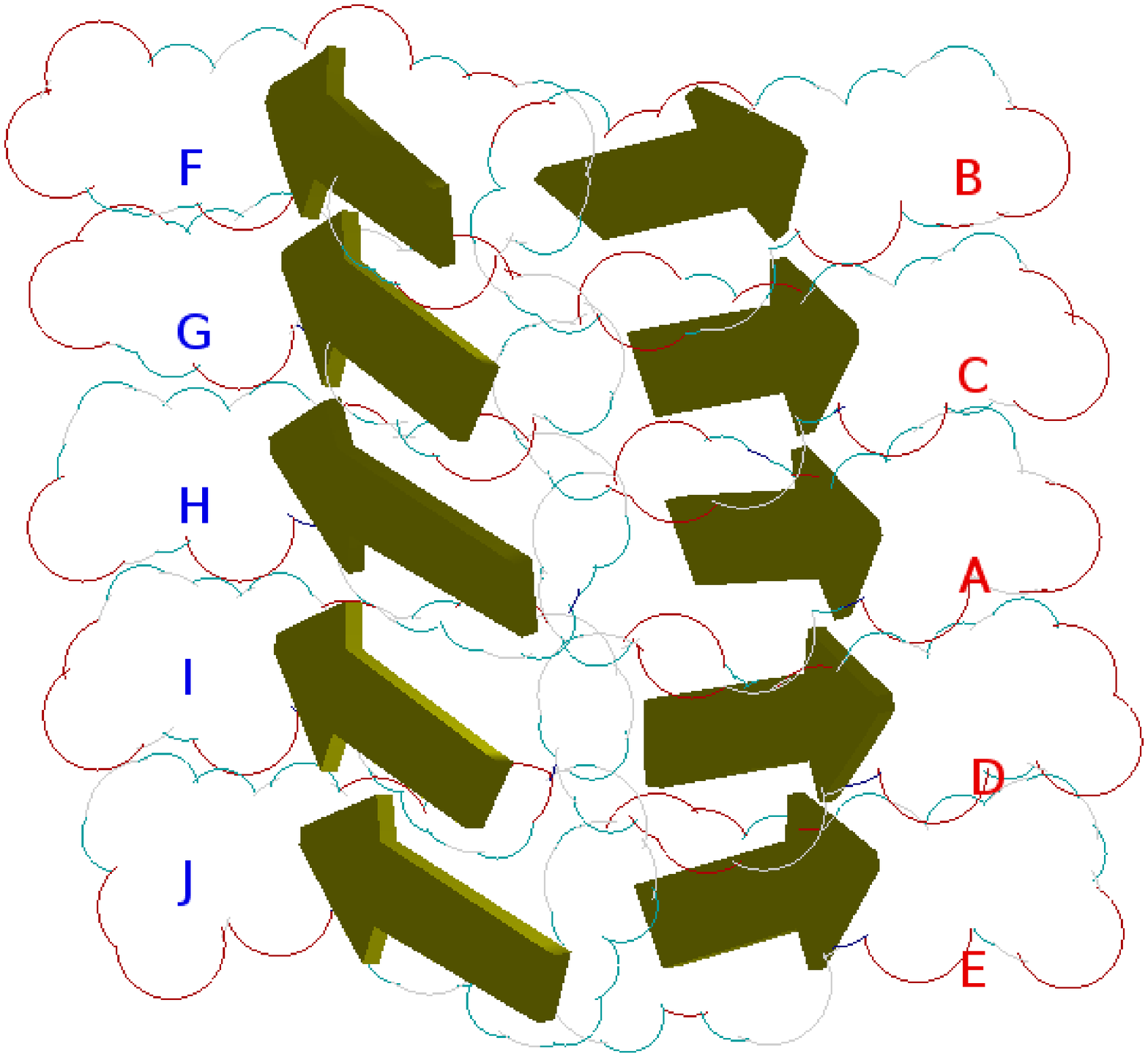}
}
\caption{Perfect 3nvh-Models 1-2 (from left to right respectively) for prion AGAAAAGA segment 113-120 . The purple dashed lines denote the hydrogen bonds. A, B, ..., I, J, K, L denote the 12 chains of the fibrils.}
\label{3nvh_CDT_models_min2}
\end{figure}

The amyloid fibril models of prion AGAAAAGA segment refined by Amber 11 are illuminated in Fig.s \ref{3nvf_CDT_models_min2}-\ref{3nvh_CDT_models_min2}. All these models are without any bad contact now (checked by package Swiss-PdbViewer), and the vdw interactions between the two $\beta$-sheets are in a very perfect way now.
All the initial structures before dealt by CDT approach have very far vdw contacts between the two $\beta$-sheets. CDT easily made the vdw contacts come closer and reach a state with the lowest potential energy, which has perfect vdw contacts as shown in Fig.s \ref{3nvf_CDT_models_min2}-\ref{3nvh_CDT_models_min2}.\\

\section{Conclusion}
Global optimization of Lennard-Jones clusters is a challenging problem for researchers in the field of biology, physics, chemistry, computer science, materials science, and especially for experts in mathematical optimization research field because of the nonconvexity of the L-J potential energy function and enormous local minima on the potential energy surface. In March 2008, American Mathematical Programming Society specially produced one whole issue, No. 76, to discuss this problem. In this paper through clever use of global optimization techniques of Gao's canonical dual theory (CDT), we successfully tackle this challenging problem illuminated by the amyloid fibril molecular model building. Clearly, this paper shows to readers that CDT is very useful and powerful to tackle challenging problems in optimization area and many other areas.\\

\noindent {\bf Acknowledgments:} This research was supported by US Air Force Office of Scientific Research under the grant AFOSR FA9550-10-1-0487, by a Victorian Life Sciences Computation Initiative (VLSCI) grant number VR0063 on its Peak Computing Facility at the University of Melbourne, an initiative of the Victorian Government, and by the Head and Colleagues of Graduate School of Ballarat University.\\

\end{document}